%% file: main.tex
\documentclass[11pt]{article}

\usepackage{graphicx}
\usepackage{latexsym,amsmath,amsfonts,amscd, amsthm, dsfont}
\usepackage{bm,color}
\usepackage{epsfig,verbatim,epstopdf,graphics}
\usepackage{subfigure}
\usepackage{changebar}
\usepackage{multirow}

\usepackage{algorithmic}

\usepackage{yhmath}
 \usepackage{booktabs} 
 \usepackage{tikz}
\usepackage{verbatim}
\usetikzlibrary{arrows,backgrounds,snakes,shapes}
 \numberwithin{equation}{section}

\graphicspath{{./}{./figure/}}
\allowdisplaybreaks

\newtheoremstyle{plainNoItalics}{}{}{\normalfont}{}{\bfseries}{.}{ }{}

\theoremstyle{plain}
\newtheorem{thm}{Theorem}[section]

\theoremstyle{plainNoItalics}

\newtheorem{exa}[thm]{Example}

\newcommand{\be}{\begin{eqnarray}}
\newcommand{\ee}{\end{eqnarray}}
\newcommand{\beno}{\begin{eqnarray*}}
\newcommand{\eeno}{\end{eqnarray*}}

\makeatletter

\newcommand{\Rmnum}[1]{\expandafter\@slowromancap\romannumeral #1@}
\makeatother

\begin{document}

\baselineskip=1.8pc


\input{title}
\input{intro}

\input{algorithm_Q}

\input{numerical}
\input{conclusion}

\appendix
\input{appendix}
\bibliographystyle{abbrv}
\bibliography{refer17}

\end{document}

%% file: title.tex
\begin{center}
{\bf
A high order semi-Lagrangian discontinuous Galerkin method for Vlasov-Poisson simulations without operator splitting
}
\end{center}

\vspace{.2in}
\centerline{
Xiaofeng Cai\footnote{
 Department of Mathematics, University of Delaware, Newark, DE, 19716. E-mail: xfcai89@gmail.com.
},
 Wei Guo \footnote{
Department of Mathematics and Statistics, Texas Tech University, Lubbock, TX, 70409. E-mail:
weimath.guo@ttu.edu. Research
is supported by NSF grant NSF-DMS-1620047.
},
Jing-Mei Qiu\footnote{Department of Mathematics, University of Delaware, Newark, DE, 19716. E-mail: jingqiu@udel.edu. Research supported by NSF grant NSF-DMS-1522777 and Air Force Office of Scientific Computing FA9550-12-0318.}
}

\bigskip
\noindent
{\bf Abstract.}
In this paper, we develop a high order semi-Lagrangian (SL) discontinuous Galerkin (DG) method for nonlinear Vlasov-Poisson (VP) simulations without operator splitting. In particular, we combine two recently developed novel techniques: one is the high order non-splitting SLDG transport method \emph{[Cai, et al.,  J Sci Comput, 2017]}, and the other is the high order characteristics tracing technique proposed in \emph{[Qiu and Russo, J Sci Comput, 2017]}. The proposed method with up to third order  accuracy in both space and time is locally mass conservative, free of splitting error, positivity-preserving, stable and robust for large time stepping size.
The SLDG VP solver is applied to classic benchmark test problems such as Landau damping and two-stream instabilities for VP simulations. Efficiency and effectiveness of the proposed scheme is extensively tested. Tremendous CPU savings are shown by comparisons between the proposed SL DG scheme and the classical Runge-Kutta DG method. 

\vfill

\newpage

%% file: intro.tex
t\section{Introduction}

This paper focuses on development of a class of high order semi-Lagrangian discontinuous Galerkin (SLDG) methods for Vlasov-Poisson (VP)
simulations without operator splitting. This is a continuation of our previous research effort on a high order non-splitting SLDG method for solving linear transport equations \cite{cai2016high}. The VP system,
arising from plasma applications, is known as a fundamental model for collisionless plasmas with a negligible magnetic field.  It reads as follows,
\begin{equation}
\frac{\partial f}{\partial t} + \mathbf{v}\cdot\nabla_{\mathbf{x}}f + \mathbf{E}(\mathbf{x},t)\cdot\nabla_{\mathbf{v}}f=0,
\label{vlasov1}
\end{equation}
and
\begin{equation}
\mathbf{E}(\mathbf{x},t) = - \nabla_{\mathbf{x}}\phi(\mathbf{x},t),\ \ -\Delta_{\mathbf{x}}\phi(\mathbf{x},t) = \rho(\mathbf{x},t),
\label{poisson}
\end{equation}
where $\mathbf{x}$ and $\mathbf{v}$ are coordinates in phase space $(\mathbf{x},\mathbf{v})\in \mathbb{R}^3\times\mathbb{R}^3$, $\mathbf{E}$ is the electric field, $\phi$ is the self-consistent electrostatic potential and $f(t,\mathbf{x},\mathbf{v})$ is probability distribution function which describes the probability of finding a particle with velocity $\mathbf{v}$ at position $\mathbf{x}$ at time $t$. The probability distribution function  couples to the long range fields via the charge density, $\rho(t,\mathbf{x}) = \int_{\mathbb{R}^3}f(t,\mathbf{x},\mathbf{v})d\mathbf{v}-1$, where we take the limit of uniformly distributed infinitely massive ions in the background. Equations \eqref{vlasov1} and \eqref{poisson} have been nondimensionalized so that all physical constants are one.

 Particle-in-cell (PIC) methods have long been a standard tool for  numerical simulation of the VP system \cite{birdsall2005plasma}. Such methods are known to be able to generate  reasonable results with relatively low computational cost for high-dimensional simulations. A PIC method mainly comprises two components: (a) A collection of $N$ macro-particles are sampled from the initial distribution function, and they are pushed in the Lagrangian framework via solving the characteristic equation
\begin{equation}
\left\{
\begin{array}{l}
\displaystyle\frac{d\mathbf{x}}{dt} = \mathbf{v},\\[2mm]
\displaystyle\frac{d\mathbf{v}}{dt} = \mathbf{E}(\mathbf{x},t).
\end{array}\right.
\label{eq:char}
\end{equation}
(b) Meanwhile, electric field $\mathbf{E}$ is solved from Poisson's equation \eqref{poisson} by means of a mesh-based scheme when needed. Despite of the simplicity of the PIC method, it suffers from the sampling noise of order $\mathcal{O}(1/\sqrt{N})$, which prevents accurate description of physics of interest in many cases. We refer to the classic textbook \cite{birdsall2005plasma} for a more detailed review of PIC methods.
In this work, we are interested in the SL approach. As with PIC methods, an SL method advances the solution by following characteristics; while, instead of particles, the solution is interpreted based on a fixed mesh,  similar to the Eularian approach. Consequently, such a method is able to
conveniently achieve desirable accuracy with the time step restriction
only set by the physical quantities such as the plasma frequency, leading to great savings
in computational time.  Due to the distinguished property, SL approaches have already elicited substantial interest in  plasma  simulation community \cite{sonnendruecker,filbet2003comparison}.

In this paper, we use the widely recognized DG spatial discretization \cite{cockburn2000discontinuous} for the VP simulation. The DG method uses a discontinuous
approximation space for the numerical solution and test functions, thus being very effective in resolving complex solution structures, such as the filaments arising from phase mixing in the Vlasov simulations.
 By contrast, a continuous  finite element method tends to introduce excessive numerical diffusion by the restrictive continuity requirement of the approximation space, resulting in smears of the solution or spurious oscillations \cite{heath2012discontinuous}. The Eularian DG methods in conjunction with the Runge-Kutta (RK) and other time integrators have  been applied for the Vlasov model in the literature, see  \cite{heath2012discontinuous,cheng2012study,cheng2014discontinuous,cheng2014energy,cheng2014energyvm}.  On the other hand, the major drawback of such Eularian DG methods is the associated stringent $CFL$ condition when an explicit time integrator is employed. The implicit method developed in \cite{cheng2014energy} can avoid the issue, but substantial computational cost would be incurred for solving the resulting nonlinear algebraic equations.

In our previous work, a class of high order SLDG methods has been proposed which incorporates  DG spatial discretization into the SL framework with the aim to take advantage of both. To our best knowledge, this is the first non-splitting, locally conservative, and highly accurate (with up to third order accuracy) SLDG scheme in the literature. In this paper, we consider to generalize  our scheme to the VP simulations. The first SL method for solving the VP system was developed by  Cheng and
Knorr in their seminal paper \cite{cheng}, in which an operator splitting strategy was introduced. 
One prominent advantage of performing operator splitting is that the resulting split equations are linear and in lower dimensions,  thus largely simplifying the algorithm design and implementation of SL schemes \cite{filbet2003comparison,casas2016high}.
Most existing high order SL schemes are designed based on the splitting strategy. In the literature, finite volume based \cite{filbet2001conservative}, finite difference based \cite{nakamura1999cubic,sonnendruecker,carrillo2007nonoscillatory,Qiu_Christlieb,qiu2011conservative}, and DG based \cite{qiu2011positivity,rossmanith2011positivity} methods are available. However, a splitting error in time will be incurred, which can be significant and hence compromise accuracy of the numerical solution over long time Vlasov simulations \cite{christlieb2014high}. This observation motivates our study on the development of a non-splitting high order SLDG schemes for the VP system. The SLDG scheme proposed in \cite{restelli2006semi} is based on a flux form and free of splitting error, but still subject to time step restriction, which undermines computational efficiency of the scheme.

When generalizing our non-splitting SLDG transport method to the VP simulation, we need to address an additional difficulty. That is, unlike in the linear transport or the VP simulation with operator splitting setting, the characteristics can no longer be exactly tracked, since their evolution is nonlinear and subject to the electric field that is induced
by the unknown function itself. In \cite{qiu2017high}, a novel characteristics tracking strategy with up to third order temporal accuracy is developed via a two-stage multi-derivative prediction-correction approach. We propose in this work to incorporate the strategy to realize high order accuracy in time. It is worthwhile to mention that, in \cite{qiu2017high}, a non-splitting finite difference SL method based on such a characteristics tracking strategy is also proposed, which is high order accurate in both space and time and unconditionally stable, but fails to conserve the total mass. More recently, a mass conservative variant is developed in \cite{xiong2016conservative} through a conservative correction technique, yet a time step constraint is introduced for stability.

This paper is organized as follows. In Section 2, we formulate the SLDG scheme for solving the VP system. We also propose an economical version of the SLDG method for computational efficiency. The performance of the proposed method is shown in Section 3 through extensive numerical tests on several benchmark problems for VP simulations. Finally, concluding remarks are made in Section 4.

%% file: algorithm_Q.tex
\section{Truly multi-dimensional SLDG algorithm }

In this section, we formulate the SLDG method for the VP simulations in two dimensions. We start by introducing the underlying algorithm framework, then the two main components of the proposed SLDG scheme are presented, including a high order characteristics tracing mechanism for the VP system based on a prediction-correction technique as well as the high order SLDG transport scheme. Some implementation details are also provided.

\subsection{Algorithm framework}
We consider the VP system \eqref{vlasov1} with the one-dimensional (1D) physical space and 1D velocity space
on the two-dimensional (2D) domain $\Omega=\Omega_x\times\Omega_v$, where $\Omega_x$ is a bounded domain with periodic boundary conditions and $\Omega_v=[-v_{\max}, v_{\max}]$ with $v_{\max}$ chosen large enough so that a zero boundary condition is reasonably imposed.
We assume a Cartesian uniform partition of the computational domain $\Omega=\{A_j\}_{j=1}^J$ (see Figure \ref{schematic_2d}) for simplicity. In principle, the method can be extended to general unstructured meshes with some modifications in implementation. Such extension will be addressed in our forthcoming paper.
We define the finite dimensional piecewise polynomial approximation space, $V_h^k = \{ v_h: v_h|_{A_j} \in P^k(A_j) \}$, where $P^k(A_j)$ denotes the set of polynomials of degree at most of $k$ on element $A_j$.
For illustrative purposes, we only present the formulation of the SLDG scheme with $P^1$ polynomial space. The generalization to $P^2$ polynomial space follows a similar procedure discussed in \cite{cai2016high}.

In order to update the solution at time level $t^{n+1}$ over the cell $A_j$ based on the solution at time level $t^n$, we employ the weak formulation of characteristic Galerkin method proposed in \cite{Guo2013discontinuous, cai2016high}. Specifically, we consider the following adjoint problem for the time dependent test function $\psi$
\begin{equation}
\psi_t + v\psi_x + E(x,t) \psi_v = 0,\
\text{subject to} \
\psi(t=t^{n+1}) = \Psi(x,v),\
t\in [t^n, t^{n+1} ],
\label{adjoint}
\end{equation}
where $\Psi\in P^k(A_j)$.
The scheme formulation takes advantage of the identity
\begin{equation}
\frac{d}{dt} \int_{  \widetilde{A}_j(t)} f(x,v,t) \psi(x,v,t) dxdv =0,
\end{equation}
where $\widetilde{A}_j(t)$ is a dynamic moving cell, emanating form the Eulerian cell $A_j$ at $t^{n+1}$ backward in time by following characteristics trajectories.
The multi-dimensional  SLDG scheme is formulated as follows:
Given the approximate solution $f^n\in V_h^k$ at time $t^n$, find $f^{n+1}\in V_h^k$ such that $\forall \Psi\in V_h^k$, we have
\begin{equation}
\int_{A_j} f^{n+1} \Psi(x,v) dxdv
=
\int_{A_j^\star} f^n \psi(x,v,t^n) dxdv, \quad \mbox{for} \quad j = 1, \cdots, J,
\label{sldg}
\end{equation}
where $\psi$ solves \eqref{adjoint} and $A_j^\star= \tilde{A}_j(t^n)$. $A_j^\star$ is called the upstream cell of $A_j$.
In general, $A_j^\star$ is no longer a rectangle, for example,
see a deformed upstream cell bounded by red curves in Figure \ref{schematic_2d}.
The proposed SLDG method in updating the numerical solution $f^n$ to $f^{n+1}$ consists of the following two main steps.
\begin{description}
  \item[1.]  {\bf Construct approximated upstream cells by following characteristics.} 
  Denote the four vertices of $A_j$  as $c_q$, with the coordinates $( x_q, v_q )$, $q=1,\cdots,4$ in the phase space.
  We trace characteristics backward in time to $t^n$ for the four vertices and then obtain $c_q^\star$ with the new coordinates $( x_q^\star, v_q^\star ),q=1,\cdots,4$. For example, see $c_4$ and $c_4^\star$ in Figure \ref{schematic_2d}.
  Then the upstream cell can be approximated by a quadrilateral determined by the four vertices $c_q^\star$.
The new coordinates $( x_q^\star, v_q^\star )$ of $c^\star_q$ are approximated by numerically solving the characteristics equation \eqref{adjoint} in the 2D case, i.e.,
      \begin{equation}
      \begin{cases}
      \frac{dx(t)}{dt} = v(t), \\[3mm]
      \frac{dv(t)}{dt} = E(x(t),t),
       \end{cases}
\quad       \mbox{with} \quad
           \begin{cases}
      x(t^{n+1}) = x_q, \\
      v(t^{n+1} ) = v_q,
     \end{cases}
     \quad
q = 1, 2, 3, 4,
      \label{characteristics}
      \end{equation}
      which is a set of final value problems.
We remark that the above equations are non-trivial to solve with high order temporal accuracy. In particular, note that the electric field $E$ depends on the unknown $f$ via Poisson's equation \eqref{poisson} in a global rather than local fashion.
  To circumvent the difficulty, we propose to combine a high order two-stage multi-derivative prediction-correction strategy for tracing characteristics as proposed in \cite{qiu2017high}. Such a strategy is briefly described in the context of the proposed SLDG scheme in Section \ref{prediction-correction}. If a high order (e.g. third order) approximation is desired, then four sides of $A_j^\star$ should be approximated by quadratic curves, in which case more points should be tracked for the curve fitting, see \cite{cai2016high} for more details.


  \item[2.]  {\bf Update the solution $f^{n+1}$ by evaluating the RHS of eq.~\eqref{sldg} for $\forall\Psi\in V_h^k$.}
We approximate $A_j^\star$ by a quadrilateral in the previous step. The test function $\psi$ at $t^n$ can be approximated by a polynomial via a least squares procedure by tracking point values of $\psi$ along characteristics.
In order to efficiently evaluate the volume integral in the RHS of \eqref{sldg}, it is converted into a set of line integrals by the use of Green's theorem. Such an idea is borrowed from CSLAM \cite{lauritzen2010conservative}, and further reformulated in \cite{cai2016high} for the development of an SLDG transport scheme. The above-mentioned procedure is briefly described in Section \ref{section:SLDG}.

\end{description}
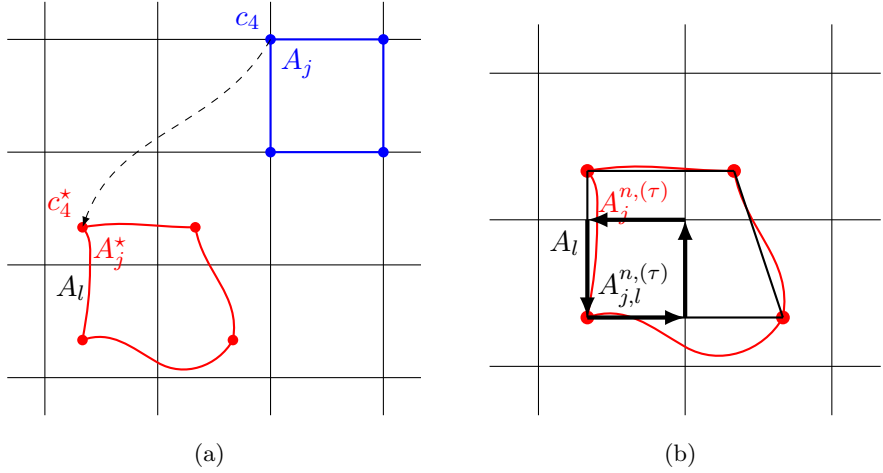
\begin{figure}[h]
\centering
\subfigure[]{
\begin{tikzpicture}
    \draw[black,thin] (0,0.5) node[left] {} -- (5.5,0.5)
                                        node[right]{};
    \draw[black,thin] (0,2.) node[left] {$$} -- (5.5,2)
                                        node[right]{};
    \draw[black,thin] (0,3.5) node[left] {$$} -- (5.5,3.5)
                                        node[right]{};
    \draw[black,thin] (0,5 ) node[left] {$$} -- (5.5,5)
                                        node[right]{};
    \draw[black,thin] (0.5,0) node[left] {} -- (0.5,5.5)
                                        node[right]{};
    \draw[black,thin] (2,0) node[left] {$$} -- (2,5.5)
                                        node[right]{};
    \draw[black,thin] (3.5,0) node[left] {$$} -- (3.5,5.5)
                                        node[right]{};
    \draw[black,thin] (5,0) node[left] {$$} -- (5,5.5)
                                        node[right]{};
    \fill [blue] (3.5,3.5) circle (2pt) node[] {};
    \fill [blue] (5,3.5) circle (2pt) node[] {};
    \fill [blue] (3.5,5) circle (2pt) node[below right] {$A_j$} node[above left] {$c_4$};
    \fill [blue] (5,5) circle (2pt) node[] {};

     \draw[thick,blue] (3.5,3.5) node[left] {} -- (3.5,5)
                                        node[right]{};
      \draw[thick,blue] (3.5,3.5) node[left] {} -- (5,3.5)
                                        node[right]{};
       \draw[thick,blue] (3.5,5) node[left] {} -- (5,5)
                                        node[right]{};
        \draw[thick,blue] (5,3.5) node[left] {} -- (5,5)
                                        node[right]{};
    \fill [red] (1.,1) circle (2pt) node[above right,black] {};
    \fill [red] (3,1) circle (2pt) node[] {};
    \fill [red] (1,2.5) circle (2pt) node[below right] {$A_j^\star$} node[above left] {$c_4^\star$};
    \fill [red] (2.5,2.5) circle (2pt) node[] {};

     \draw[-latex,dashed](3.5,5)node[right,scale=1.0]{}
        to[out=240,in=70] (1,2.50) node[] {};

     \draw (0.5+0.01,2-0.01) node[fill=white,below right] {$A_l$};

     \draw [red,thick] (1,1)node[right,scale=1.0]{}
        to[out=20,in=150] (2,0.7) node[] {};

        \draw [red,thick] (2,0.7)node[right,scale=1.0]{}
        to[out=330,in=240] (3,1) node[] {};
             \draw [red,thick] (1,2.5)node[right,scale=1.0]{}
        to[out=310,in=90] (1.1,2) node[] {};
        \draw [red,thick] (1.1,2)node[right,scale=1.0]{}
        to[out=270,in=80] (1,1) node[] {};

        \draw [red,thick] (1,2.5)node[right,scale=1.0]{}
        to[out=10,in=180] (2.5,2.5) node[] {};

        \draw [red,thick] (3,1)node[right,scale=1.0]{}
        to[out=80,in=280] (2.5,2.5) node[] {};
\end{tikzpicture}
}
\subfigure[]{

\begin{tikzpicture}[scale = 1.3]
    \draw[black,thin] (0,0.5) node[left] {} -- (4,0.5)
                                        node[right]{};
    \draw[black,thin] (0,2.) node[left] {$$} -- (4,2)
                                        node[right]{};
    \draw[black,thin] (0,3.5) node[left] {$$} -- (4,3.5)
                                        node[right]{};
    \draw[black,thin] (0.5,0) node[left] {} -- (0.5,4)
                                        node[right]{};
    \draw[black,thin] (2,0) node[left] {$$} -- (2,4)
                                        node[right]{};
    \draw[black,thin] (3.5,0) node[left] {$$} -- (3.5,4)
                                        node[right]{};

    \fill [red] (1.,1) circle (2pt) node[above right,black] {$A_{j,l}^{n,(\tau)}$};
    \fill [red] (3,1) circle (2pt) node[] {};
    \fill [red] (1,2.5) circle (2pt) node[below right] {$A_j^{n,(\tau)}$} node[above left] {};
    \fill [red] (2.5,2.5) circle (2pt) node[] {};

     \draw (0.5+0.01,2-0.01) node[fill=white,below right] {};

\draw [red,thick] (1,1)node[right,scale=1.0]{}
        to[out=20,in=150] (2,0.7) node[] {};

        \draw [red,thick] (2,0.7)node[right,scale=1.0]{}
        to[out=330,in=240] (3,1) node[] {};
             \draw [red,thick] (1,2.5)node[right,scale=1.0]{}
        to[out=310,in=90] (1.1,2) node[] {};
        \draw [red,thick] (1.1,2)node[right,scale=1.0]{}
        to[out=270,in=80] (1,1) node[] {};

        \draw [red,thick] (1,2.5)node[right,scale=1.0]{}
        to[out=10,in=180] (2.5,2.5) node[] {};

        \draw [red,thick] (3,1)node[right,scale=1.0]{}
        to[out=80,in=280] (2.5,2.5) node[] {};
           \draw (0.5+0.01,2-0.01) node[fill=white,below right] {$A_l$};
         \draw[-latex,ultra thick] (1,1)node[right,scale=1.0]{}
        to  (2,1) node[] {};

             \draw[-latex,ultra thick]  (2,1)node[right,scale=1.0]{}
        to (2,2) node[] {};
             \draw[-latex,ultra thick]  (2,2)node[right,scale=1.0]{}
        to (1,2) node[] {};
             \draw[-latex,ultra thick]   (1,2)node[right,scale=1.0]{}
        to (1,1) node[] {};

\draw [thick] (1,1)-- (3,1) node[] {};
 \draw [thick] (1,2.5) -- (1,1) node[] {};
        \draw [thick] (1,2.5)--(2.5,2.5) node[] {};
        \draw [thick] (3,1)--(2.5,2.5) node[] {};
\end{tikzpicture}

}
\caption{Schematic illustration of the SLDG formulation in two dimension: quadrilateral approximation to a upstream cell.  }
\label{schematic_2d}
\end{figure}

\subsection{High order characteristics tracing prediction-correction algorithm}
\label{prediction-correction}

In this subsection,  we describe a high order predictor-corrector procedure for locating the feet of the characteristics of the VP system. Such an approach is originally proposed in \cite{qiu2017high}.
We first introduce several shorthand notations. The superscript $^n$ denotes the time level, the superscript $^{(\tau)}$ denotes the formal order of approximation for time discretization, and the subscript $q$ is the index for the vertices of the underlying cell  in the phase space. For example, $(x_q^{n, (\tau)}, v_q^{n, (\tau)})$ is the $\tau$-th order approximation of $(x_q^\star, v_q^\star)$ and $A_j^{n, (\tau)}$ is the quadrilateral determined by the corresponding four vertices.

We start from a first order scheme for tracing characteristics \eqref{characteristics}. We let
\begin{equation}
x_{q }^{n,(1)}
=
x_{q} - v_{q}\Delta t, \
v_{q}^{n,(1)}
=
v_{q} - E(x_{q},t^n) \Delta t, \
\label{eq:first}
\end{equation}
which leads to a first order approximations to $(x_{q}^\star,v_{q}^\star)$.
The electric field $E$ depends on $f$ via Poisson's equation, which can be numerically solved by the LDG method \cite{arnold2002unified, cockburn1998local, castillo2000priori}.
Note that the numerical solution $E_h$ solved by the LDG method are discontinuous across cell boundaries, i.e. $E_h(x_q^+,t^n)\neq E_h(x_q^-,t^n)$, where the superscripts $^\pm$ are the right- and left-hand limits of the corresponding functions with respect to $x_q$. In our implementation, we take the average of $E_h$ at the boundaries as the numerical solution $E_q$, i.e.
$E(x_{q},t^n) = \frac{(E_h(x_q^+,t^n) + E_h(x_q^-,t^n) )}{2}$.
Let $A_j^{n,(1)}$ to be the quadrilateral formed by the four upstream vertices $(x_{q }^{n,(1)}, v_{q }^{n,(1)})$, $q=1, 2, 3, 4.$
Then, by the SLDG formulation (to be described in the next subsection)
 \begin{equation}
   \int_{A_j} f^{n+1,(1)} \Psi(x,v) dxdv = \int_{ A_j^{n,(1)} } f^n \psi(x,v,t^n) dxdv,
   \label{eq: SLDG_t1}
  \end{equation}
we obtain $f^{n+1, (1)}$ as a first order approximation in time to the solution at $t^{n+1}$.
Based on $f^{n+1, (1)}$, we apply the LDG method to Poisson's equation \eqref{poisson} again and compute $E_q^{n+1,(1)}$, which approximates $E(x_q, t^{n+1})$ with first order temporal accuracy.

A second order scheme can be built upon the first order one. First,  let
\begin{equation}
x_{q }^{n,(2)}
=
x_{q} - \frac12\left( v_{q}  + v_{q}^{n,(1)}  \right)\Delta t, \
v_{q}^{n,(2)}
=
v_{q} - \frac12 \left(  E( x_{q }^{n,(1)},t^n )  + E_{q}^{n+1,(1)}   \right) \Delta t,\
\label{eq:second}
\end{equation}
which gives a second order approximations to  $(x_q^\star,v_q^\star)$.
Then the second order approximation solution $f^{n+1,(2)}$ is obtained from the SLDG formulation
 \begin{equation}
   \int_{A_j} f^{n+1,(2)} \Psi(x,v) dxdv = \int_{ A_j^{n,(2)} } f^n \psi(x,v,t^n) dxdv.
   \label{eq: SLDG_t2}
  \end{equation}
Based on $f^{n+1,(2)}$, we are able to compute $E_q^{n+1,(2)}$ from Poisson's equation, which approximates $E(x_q,t^{n+1})$ with second order temporal accuracy.

Lastly, a third order scheme can be designed based on the above second order approximation. Let
\begin{align}
x_{q }^{n,(3)}
=
x_{q} -
v_{q} \Delta t
+\frac{\Delta t^2}{2}
\left( \frac23 E_{q}^{n+1,(2)} + \frac13 E( x_{q}^{n,(2)} ,t^n  )  \right),
\label{eq:third1}\\
v_{q}^{n,(3)}
=
v_{q}
-
E^{n+1,(2)}_{q} \Delta t
+
\frac{\Delta t^2}{2}
\left(
\frac23( \frac{d}{dt}E(x_{q} ,t^{n+1} ) )^{(2)} + \frac13 \frac{d}{dt} E( x_{q}^{n,(2)},t^n )
\right),
\label{eq:third2}
\end{align}
where  $\frac{d}{dt}$ is the material derivative along the characteristic curve \cite{qiu2017high}, i.e.,
$$\frac{d}{dt} E=\frac{\partial E}{\partial t} +\frac{\partial E}{\partial x}\frac{\partial x}{\partial t} = \overline{J}^0 - J(x,t) + v(\rho-1).$$
Note that we have used Amp\`{e}re's law
$$\frac{\partial E}{\partial t} = \overline{J}^0 - J(x,t),$$
where $J(x,t) = \int fvdv$ is the current density and $\overline{J}^0 = \frac{1}{L} \int J(x,t=0)dx$ is the average of the current density at $t=0$ with $L=|\Omega_x|$.
In particular, in equation \eqref{eq:third2}
\begin{align}
\left(  \frac{d}{dt} E(x_q,t^{n+1}) \right)^{(2)}
& =
\overline{J}^0  - J_q^{n+1,(2)} + v_q( \rho_q^{n+1,(2)} - 1 ),\\
\frac{d}{dt} E( x_q^{n,(2)} , t^n )
& =
\overline{J}^0 - J( x_q^{n,(2)} )  + v_q^{n,(2)} (  \rho( x_q^{n,(2)} ,t^n ) -1 ).
\end{align}
It can be checked by a local truncation error analysis that $(x_{q }^{n,(3)}, v_{q }^{n,(3)})$ is a third order approximation to $(x_q^\star,v_q^\star)$ \cite{qiu2017high}. Consequently, the third order approximation solution
$f^{n+1,(2)}$ is updated from the SLDG formulation
\begin{equation}
\int_{A_j} f^{n+1,(2)} \Psi(x,v) dxdv = \int_{ A_j^{n,(3)} } f^n \psi(x,v,t^n) dxdv.
\label{eq: SLDG_t3}
\end{equation}

%
%

\subsection{A two-dimensional SLDG method with quadrilateral upstream cells.}
\label{section:SLDG}

Below, we present the procedure in evaluating the integral $\int_{A_j^\star} f^n \psi(x, v, t^n) dxdv$ with a quadrilateral upstream cell $A_j^\star$. In the algorithm design, we have to pay attention to the following two observations, see \cite{cai2016high}.
\begin{itemize}
\item $\Psi = \psi(x,v, t^{n+1})$ is chosen to be polynomial basis functions on $V_h^k$, while, in general $\psi(x, v, t^n)$ is no longer a polynomial. A polynomial function constructed by a least squares procedure is used to approximate $\psi(x, v, t^n)$.
\item  Over the upstream cell $A_j^\star$ (or its approximation $A_{j}^{n,(\tau)}$), $f^n(x, v, t^n)$ is discontinuous across Eulerian cell boundaries, see the background Eulerian grid lines in Figure~\ref{schematic_2d}. To properly evaluate the volume integral, one has to perform the evaluation in a sub-area by sub-area manner. Meanwhile, direct evaluation of volume integrals over these irregular-shape sub-areas is very involved in implementation. The proposed strategy is to convert each volume integral into line integrals by the use of Green's Theorem.
\end{itemize}
Based on these observations, the proposed algorithm consists of two main components. One is the search algorithm that finds the boundaries for each sub-area, i.e. the overlapping region between the upstream cell and background Eulerian cells. The other is the use of Green's theorem that enables us to convert the volume integral to line integrals based on the result of the search algorithm. Below we describe the detailed procedure in evaluating the volume integral over an approximation of upstream cell $A_{j}^{n,(\tau)}$ for the SLDG scheme with $P^1$ polynomial space.

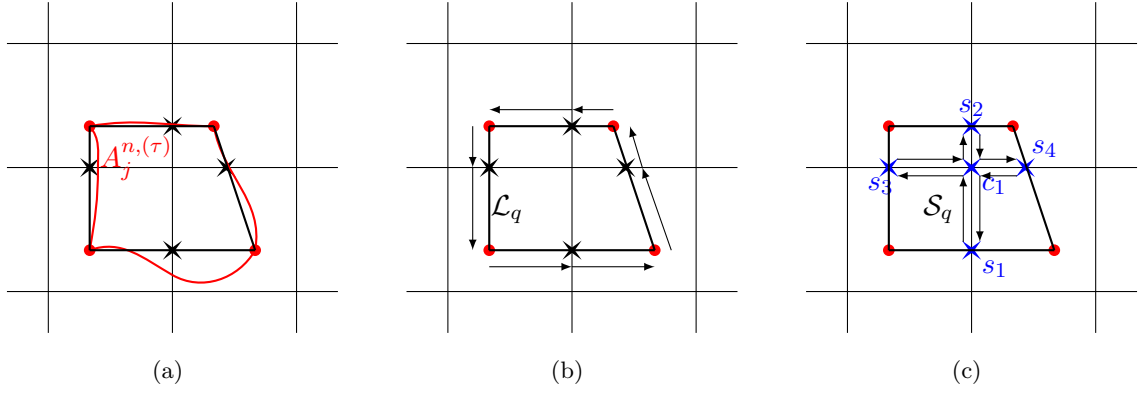
\begin{figure}[h]
\centering

\subfigure[]{
\begin{tikzpicture}[scale = 1.1]
    \draw[black,thin] (0,0.5) node[left] {} -- (4,0.5)
                                        node[right]{};
    \draw[black,thin] (0,2.) node[left] {$$} -- (4,2)
                                        node[right]{};
    \draw[black,thin] (0,3.5) node[left] {$$} -- (4,3.5)
                                        node[right]{};
    \draw[black,thin] (0.5,0) node[left] {} -- (0.5,4)
                                        node[right]{};
    \draw[black,thin] (2,0) node[left] {$$} -- (2,4)
                                        node[right]{};
    \draw[black,thin] (3.5,0) node[left] {$$} -- (3.5,4)
                                        node[right]{};

    \fill [red] (1.,1) circle (2pt) node[above right,black] {};
    \fill [red] (3,1) circle (2pt) node[] {};
    \fill [red] (1,2.5) circle (2pt) node[below right] {$A_j^{n,(\tau)}$} node[above left] {};
    \fill [red] (2.5,2.5) circle (2pt) node[] {};
   \usetikzlibrary{shapes.geometric}
  \node[fill,star,star points=4, star point ratio=.2] at (2,1) {};
  \node[fill,star,star points=4, star point ratio=.2] at (2,2.5) {};
  \node[fill,star,star points=4, star point ratio=.2] at (1,2) {};
  \node[fill,star,star points=4, star point ratio=.2] at (2.65,2) {};

     \draw (0.5+0.01,2-0.01) node[fill=white,below right] {};
\draw [red,thick] (1,1)node[right,scale=1.0]{}
        to[out=20,in=150] (2,0.7) node[] {};

        \draw [red,thick] (2,0.7)node[right,scale=1.0]{}
        to[out=330,in=240] (3,1) node[] {};
             \draw [red,thick] (1,2.5)node[right,scale=1.0]{}
        to[out=310,in=90] (1.1,2) node[] {};
        \draw [red,thick] (1.1,2)node[right,scale=1.0]{}
        to[out=270,in=80] (1,1) node[] {};

        \draw [red,thick] (1,2.5)node[right,scale=1.0]{}
        to[out=10,in=180] (2.5,2.5) node[] {};

        \draw [red,thick] (3,1)node[right,scale=1.0]{}
        to[out=80,in=280] (2.5,2.5) node[] {};
\draw [thick] (1,1)-- (3,1) node[] {};
 \draw [thick] (1,2.5) -- (1,1) node[] {};
        \draw [thick] (1,2.5)--(2.5,2.5) node[] {};
        \draw [thick] (3,1)--(2.5,2.5) node[] {};

\end{tikzpicture}

}
\subfigure[]{
\begin{tikzpicture}[scale = 1.1]
    \draw[black,thin] (0,0.5) node[left] {} -- (4,0.5)
                                        node[right]{};
    \draw[black,thin] (0,2.) node[left] {$$} -- (4,2)
                                        node[right]{};
    \draw[black,thin] (0,3.5) node[left] {$$} -- (4,3.5)
                                        node[right]{};
    \draw[black,thin] (0.5,0) node[left] {} -- (0.5,4)
                                        node[right]{};
    \draw[black,thin] (2,0) node[left] {$$} -- (2,4)
                                        node[right]{};
    \draw[black,thin] (3.5,0) node[left] {$$} -- (3.5,4)
                                        node[right]{};

    \fill [red] (1.,1) circle (2pt) node[above right,black] {};
    \fill [red] (3,1) circle (2pt) node[] {};
    \fill [red] (1,2.5) circle (2pt) node[below right] {} node[above left] {};
    \fill [red] (2.5,2.5) circle (2pt) node[] {};

     \draw (0.5+0.01,2-0.01) node[fill=white,below right] {};

\draw [thick] (1,1)-- (3,1) node[] {};
 \draw [thick] (1,2.5) -- (1,1) node[] {};
        \draw [thick] (1,2.5)--(2.5,2.5) node[] {};
        \draw [thick] (3,1)--(2.5,2.5) node[] {};
   \usetikzlibrary{shapes.geometric}
  \node[fill,star,star points=4, star point ratio=.2] at (2,1) {};
  \node[fill,star,star points=4, star point ratio=.2] at (2,2.5) {};
  \node[fill,star,star points=4, star point ratio=.2] at (1,2) {};
  \node[fill,star,star points=4, star point ratio=.2] at (2.65,2) {};
  \draw [-latex] (1-0.2,2) -- node[right=3pt]{$\mathcal{L}_q$}(1-0.2,1) node[] {};
  \draw [-latex] (1-0.2,2.5) -- (1-0.2,2) node[] {};

  \draw [-latex] (2.5,2.5+0.2) -- (2,2.5+0.2);
  \draw [-latex] (2,2.5+0.2) -- (1,2.5+0.2);

  \draw [-latex] (1,1-0.2) -- (2,1-0.2) node[] {};
  \draw [-latex] (2,1-0.2) -- (3,1-0.2) node[] {};
  \draw [-latex] (3+0.2,1) -- (2.65+0.2,2) node[] {};
  \draw [-latex] (2.65+0.2,2) -- (2.5+0.2,2.5) node[] {};



\end{tikzpicture}

}
\subfigure[]{
\begin{tikzpicture}[scale = 1.1]
%

\node [below right,blue] at (2,1) {$s_1$};
\node [above,blue] at (2,2.5) {$s_2$};
\node [below,blue] at (0.9,2.) {$s_3$};
\node [above right,blue] at (2.6,2) {$s_4$};
\node [below right, blue] at (2,2) {$c_1$};
    \draw[black,thin] (0,0.5) node[left] {} -- (4,0.5)
                                        node[right]{};
    \draw[black,thin] (0,2.) node[left] {$$} -- (4,2)
                                        node[right]{};
    \draw[black,thin] (0,3.5) node[left] {$$} -- (4,3.5)
                                        node[right]{};
    \draw[black,thin] (0.5,0) node[left] {} -- (0.5,4)
                                        node[right]{};
    \draw[black,thin] (2,0) node[left] {$$} -- (2,4)
                                        node[right]{};
    \draw[black,thin] (3.5,0) node[left] {$$} -- (3.5,4)
                                        node[right]{};

    \fill [red] (1.,1) circle (2pt) node[above right,black] {};
    \fill [red] (3,1) circle (2pt) node[] {};
    \fill [red] (1,2.5) circle (2pt) node[below right] {} node[above left] {};
    \fill [red] (2.5,2.5) circle (2pt) node[] {};

   \usetikzlibrary{shapes.geometric}
  \node[fill,star,star points=4, star point ratio=.2,blue] at (2,1) {};
  \node[fill,star,star points=4, star point ratio=.2,blue] at (2,2.5) {};
  \node[fill,star,star points=4, star point ratio=.2,blue] at (1,2) {};
  \node[fill,star,star points=4, star point ratio=.2,blue] at (2.65,2) {};
  \node[fill,star,star points=4, star point ratio=.2,blue] at (2,2) {};

  \draw [-latex] (2-0.1,2+0.1) -- (2-0.1,2.5-0.1) node[] {};
  \draw [-latex] (2+0.1,2.5-0.1)--(2+0.1,2+0.1 )  node[] {};

   \draw [-latex] (2-0.1,1+0.1) --node[auto]{$\mathcal{S}_q$} (2-0.1,2-0.1) node[] {};
  \draw [-latex] (2+0.1,2-0.1 )--(2+0.1,1+0.1 )  node[] {};
  \draw [-latex] (2.65-0.1,2-0.1) --(2+0.1,2-0.1)  node[] {};
  \draw [-latex] (2+0.1,2+0.1)--(2.65-0.1,2+0.1)  node[] {};
    \draw [-latex] (1+0.1,2+0.1) --(2-0.1,2+0.1)  node[] {};
  \draw [-latex] (2-0.1,2-0.1)--(1+0.1,2-0.1)  node[] {};

     \draw (0.5+0.01,2-0.01) node[fill=white,below right] {};

\draw [thick] (1,1)-- (3,1) node[] {};
 \draw [thick] (1,2.5) -- (1,1) node[] {};
        \draw [thick] (1,2.5)--(2.5,2.5) node[] {};
        \draw [thick] (3,1)--(2.5,2.5) node[] {};
\end{tikzpicture}

}
\caption{Schematic illustration of the search algorithm. }
\label{schematic_search}
\end{figure}

\begin{description}
   \item[(1)] \textbf{\emph{Least squares approximation of test function $\psi(x, v, t^n)$.}}
              Based on the fact that the solution of the adjoint problem \eqref{adjoint} stays unchanged along characteristics, we have
              \begin{equation*}
              \psi( x_q^{n,(\tau)}  ,  v_q^{n,(\tau)}   , t^n ) = \Psi( x_q, v_q ),\ \ q=1,2,\cdots,4.
              \end{equation*}
              Thus, we can reconstruct a unique linear function $\psi^\star(x,v)$ by a least squares strategy that approximates $\psi(x,v,t^n)$ on $A_{j}^{n,(\tau)}$ .

   \item[(2)] 
   \textbf{\emph{Evaluation of the volume integral.}}
   Denote $A_{j,l}^{n,(\tau)}$ as a non-empty overlapping region between the upstream cell $A_j^{n,(\tau)}$ and the background Eulerian cell $A_l$, i.e., $A_{j,l}^{n,(\tau)} = A_{j}^{n,(\tau)} \cap A_l$,
   see Figure \ref{schematic_2d} (b). Then the volume integral, e.g. RHS of eq.~\eqref{eq: SLDG_t1} with $\tau=1$,  becomes
      \begin{equation}
             \int_{A_j^{n,(\tau)}} f(x,v,t^{n} )\psi(x,v,t^{n} ) dxdv
           \approx
          \sum_{l\in \varepsilon_j^{n,(\tau)}}^{  } \int_{A_{j,l}^{n,(\tau)} } f(x,v,t^n)\psi^\star(x,v)dxdv,
       \label{temp1}
      \end{equation}
      where $\varepsilon_j^{n,(\tau)} =\{ l| A_{j,l}^{n,(\tau)}  \neq \emptyset \}$.
   Note that the integrands on the RHS of \eqref{temp1} are piecewise   polynomials.
       By introducing two auxiliary function $P(x,v)$ and $Q(x,v)$ such that
       \begin{equation*}
       -\frac{\partial P }{\partial v } + \frac{\partial Q}{\partial x }  =  f(x,v,t^n)\psi^\star(x,v),
       \end{equation*}
   the area integral $ \int_{A_{j,l}^{n,(\tau)} } f(x,v,t^n)\psi^\star(x,v)dxdv  $ can be converted into line integrals via Green's theorem, i.e.,
   \begin{equation}
      \int_{A_{j,l}^{n,(\tau)} } f(x,v,t^n)\psi^\star(x,v)dxdv = \oint_{\partial A_{j,l}^{n,(\tau)}}  Pdx + Qdv,
      \label{Green}
   \end{equation}
   see Figure \ref{schematic_2d} (b). Note that the choices of $P$ and $Q$ are not unique, but the value of the line integrals is independent of the choices. In the implementation, we follow the same procedure in \cite{lauritzen2010conservative} when choosing $P$ and $Q$.
    In summary, combining \eqref{temp1} and \eqref{Green}, we have the following
\begin{align}
\int_{A_{j}^{n,(\tau)}   } f(x,v,t_{n} )\psi(x,v,t_{n} ) dxdv
=&
\sum_{l\in \varepsilon_j^{n,(\tau)}  }^{  }\int_{A_{j,l}^{n,(\tau)} } u(x,v,t_{n} )\psi^\star(x,v ) dxdv   \notag \\
=&
\sum_{l\in \varepsilon_j^{n,(\tau)} }^{  } \oint_{\partial A_{j,l}^{n,(\tau)}   }  Pdx + Qdv \notag \\
=&
\sum_{q=1}^{N_o}
\int_{ \mathcal{L}_q } [P  dx +Q  d v ]  + \sum_{q=1}^{N_i}
\int_{ \mathcal{S}_q } [Pdx +Q d v ].
\label{line}
\end{align}
Note that in the above computation, we have organized the liner integrals into two categories: along outer line segments (see Figure \ref{schematic_search} (b)) and along inner line segments (see Figure \ref{schematic_search} (c)). Line segments can be uniquely determined by two end points, which are intersection points of the four sides of the upstream cell with grid lines.
We compute all intersection points and connect them in a counterclockwise orientation to obtain outer line segments, denoted as $\mathcal{L}_q$, $q=1,\cdots,N_o$, see Figure \ref{schematic_search} (b).
The line segments that are aligned with grid lines and enclosed by $A_j^{n,(\tau)}$ are defined as inner line segments, see Figure \ref{schematic_search} (c). Note that there are two orientations along each inner segment, but the corresponding line integrals have to be evaluated in their own sub-area, given that $f^n$ is discontinuous across a inner line segment. For instance, $\overrightarrow{s_1c_1}$ belongs to the left background cell and $\overrightarrow{c_1s_1}$ belongs to the right background cell. Again we refer to \cite{lauritzen2010conservative, cai2016high} for more details in implementation and in generalization to the SLDG scheme with $P^2$ polynomial space.
 \end{description}

\subsection{A two-dimensional SLDG method with quadratic-curved quadrilateral upstream cells.}

Note that the aforementioned SLDG methods with quadrilateral approximation yields the second order accuracy even though $P^2$ approximation space is used.
 In order to achieve a formal third order accuracy, a quadratic-curved quadrilateral is used to approximate each upstream cell $A^{\star}$ when  evaluating of the integral $\int_{A^{\star} } f^n \psi(x,v,t^n) dxdv$.
In particular, one can construct a parabola to approximate each side of an upstream cell.
Since the procedure of the SLDG scheme with quadratic-curved quadrilateral upstream cells is similar to that of the SLDG scheme with quadrilateral upstream cells, we only highlight the evaluation of line integrals along outer line segments.
The procedure of evaluation consists of the following main steps.

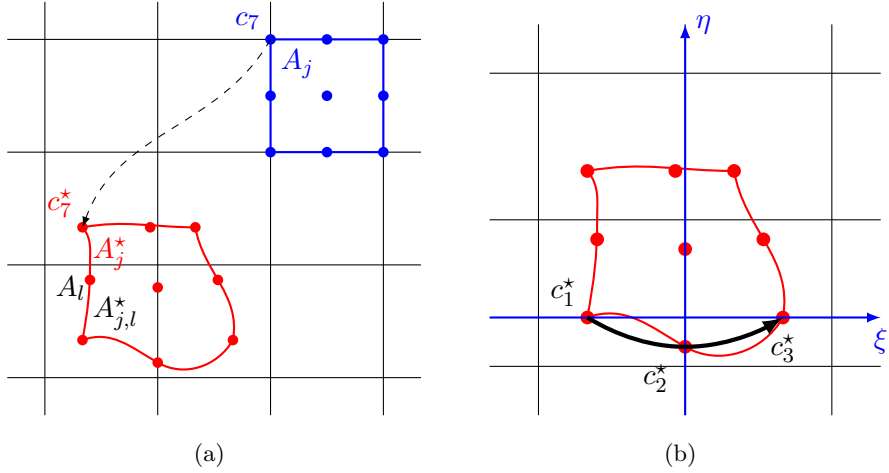
\begin{figure}[h]
\centering
\subfigure[]{
\begin{tikzpicture}
    \draw[black,thin] (0,0.5) node[left] {} -- (5.5,0.5)
                                        node[right]{};
    \draw[black,thin] (0,2.) node[left] {$$} -- (5.5,2)
                                        node[right]{};
    \draw[black,thin] (0,3.5) node[left] {$$} -- (5.5,3.5)
                                        node[right]{};
    \draw[black,thin] (0,5 ) node[left] {$$} -- (5.5,5)
                                        node[right]{};
    \draw[black,thin] (0.5,0) node[left] {} -- (0.5,5.5)
                                        node[right]{};
    \draw[black,thin] (2,0) node[left] {$$} -- (2,5.5)
                                        node[right]{};
    \draw[black,thin] (3.5,0) node[left] {$$} -- (3.5,5.5)
                                        node[right]{};
    \draw[black,thin] (5,0) node[left] {$$} -- (5,5.5)
                                        node[right]{};
    \fill [blue] (3.5,3.5) circle (2pt) node[] {};
    \fill [blue] (5,3.5) circle (2pt) node[] {};
    \fill [blue] (3.5,5) circle (2pt) node[below right] {$A_j$} node[above left] {$c_7$};
    \fill [blue] (5,5) circle (2pt) node[] {};

    \fill [blue] (3.5,4.25) circle (2pt) node[] {};
    \fill [blue] (5,4.25) circle (2pt) node[] {};
    \fill [blue] (4.25,4.25) circle (2pt) node[] {};
        \fill [blue] (4.25, 3.5) circle (2pt) node[] {};
    \fill [blue] (4.25,5) circle (2pt) node[] {};

     \draw[thick,blue] (3.5,3.5) node[left] {} -- (3.5,5)
                                        node[right]{};
      \draw[thick,blue] (3.5,3.5) node[left] {} -- (5,3.5)
                                        node[right]{};
       \draw[thick,blue] (3.5,5) node[left] {} -- (5,5)
                                        node[right]{};
        \draw[thick,blue] (5,3.5) node[left] {} -- (5,5)
                                        node[right]{};
    \fill [red] (1.,1) circle (2pt) node[above right,black] {$A_{j,l}^\star$};
    \fill [red] (3,1) circle (2pt) node[] {};
    \fill [red] (1,2.5) circle (2pt) node[below right] {$A_j^\star$} node[above left] {$c_7^\star$};
    \fill [red] (2.5,2.5) circle (2pt) node[] {};

     \draw[-latex,dashed](3.5,5)node[right,scale=1.0]{}
        to[out=240,in=70] (1,2.50) node[] {};

     \draw (0.5+0.01,2-0.01) node[fill=white,below right] {$A_l$};

     \draw [red,thick] (1,1)node[right,scale=1.0]{}
        to[out=20,in=150] (2,0.7) node[] {};

        \draw [red,thick] (2,0.7)node[right,scale=1.0]{}
        to[out=330,in=240] (3,1) node[] {};
             \draw [red,thick] (1,2.5)node[right,scale=1.0]{}
        to[out=310,in=90] (1.1,2) node[] {};
        \draw [red,thick] (1.1,2)node[right,scale=1.0]{}
        to[out=270,in=80] (1,1) node[] {};

        \draw [red,thick] (1,2.5)node[right,scale=1.0]{}
        to[out=10,in=180] (2.5,2.5) node[] {};

        \draw [red,thick] (3,1)node[right,scale=1.0]{}
        to[out=80,in=280] (2.5,2.5) node[] {};

       \fill [red] (2,0.7) circle (2pt) node[above right,black] {};
    \fill [red] (2,1.7) circle (2pt) node[] {};
    \fill [red] (1.9,2.5) circle (2pt) node[below right] {} node[above left] {};
    \fill [red] (1.1,1.8) circle (2pt) node[] {};
\fill [red] (2.8,1.8) circle (2pt) node[] {};

\end{tikzpicture}
}
\subfigure[]{

\begin{tikzpicture}[scale = 1.3]
    \draw[black,thin] (0,0.5) node[left] {} -- (4,0.5)
                                        node[right]{};
    \draw[black,thin] (0,2.) node[left] {$$} -- (4,2)
                                        node[right]{};
    \draw[black,thin] (0,3.5) node[left] {$$} -- (4,3.5)
                                        node[right]{};
    \draw[black,thin] (0.5,0) node[left] {} -- (0.5,4)
                                        node[right]{};
    \draw[black,thin] (2,0) node[left] {$$} -- (2,4)
                                        node[right]{};
    \draw[black,thin] (3.5,0) node[left] {$$} -- (3.5,4)
                                        node[right]{};

    \fill [red] (1.,1) circle (2pt) node[above right,black] {};
    \fill [red] (3,1) circle (2pt) node[] {};
    \fill [red] (1,2.5) circle (2pt) node[below right] {} node[above left] {};
    \fill [red] (2.5,2.5) circle (2pt) node[] {};

     \draw (0.5+0.01,2-0.01) node[fill=white,below right] {};

\draw [red,thick] (1,1)node[right,scale=1.0]{}
        to[out=20,in=150] (2,0.7) node[] {};

        \draw [red,thick] (2,0.7)node[right,scale=1.0]{}
        to[out=330,in=240] (3,1) node[] {};
             \draw [red,thick] (1,2.5)node[right,scale=1.0]{}
        to[out=310,in=90] (1.1,2) node[] {};
        \draw [red,thick] (1.1,2)node[right,scale=1.0]{}
        to[out=270,in=80] (1,1) node[] {};

        \draw [red,thick] (1,2.5)node[right,scale=1.0]{}
        to[out=10,in=180] (2.5,2.5) node[] {};

        \draw [red,thick] (3,1)node[right,scale=1.0]{}
        to[out=80,in=280] (2.5,2.5) node[] {};

               \fill [red] (2,0.7) circle (2pt) node[above right,black] {};
    \fill [red] (2,1.7) circle (2pt) node[] {};
    \fill [red] (1.9,2.5) circle (2pt) node[below right] {} node[above left] {};
    \fill [red] (1.1,1.8) circle (2pt) node[] {};
\fill [red] (2.8,1.8) circle (2pt) node[] {};

\draw[-latex,blue,thick](0,1)node[right,scale=1.0]{}
        to (4,1) node[below] {$\xi$};
\draw[-latex,blue,thick](2,0)node[right,scale=1.0]{}
        to (2,4) node[right] {$\eta$};

       \draw[black, ultra thick] (2,0.7) node[below left =2pt] {$c_2^\star$}
        parabola(1,1)node[above left,scale=1.0]{ $c_1^\star$ };
   \draw[-latex,black, ultra thick](2,0.7)node[above left,scale=1.0]{  }
         parabola (3,1) node[below =2pt] {$c_3^\star$};
\end{tikzpicture}

}
\caption{Schematic illustration of the SLDG formulation in two dimensions: quadratic-curved quadrilateral approximation to an upstream cell.}
\label{schematic_2d_p2}
\end{figure}

 \begin{description}
  \item[(1)] \textbf{\emph{Construction of quadratic-curved upstream cells by following characteristics.}}
   \begin{description}
     \item[(1a)]Locate nine points on the upstream cell $A_{j}^\star$, i.e.  $c_i^\star,\,i=1,\ldots,9$, from nine uniformly distributed points at $A_j$ by tracking characteristics backward in time, i.e., solving the characteristics equations \eqref{characteristics} starting from $c_i$, $i=1 \cdots, 9$, see the layout in Figure \ref{schematic_2d_p2} (a).
   \item[(1b)]  Construct a quadratic curve to approximate each side of the upstream cell. In particular, to construct the quadratic curve, $\wideparen{ c_1^\star,c_2^\star, c_3^\star }$ as shown in Figure \ref{schematic_2d_p2} (b), we perform the following procedure.
     First, we can construct  a linear coordinate transformation $x-v$ to $\xi-\eta$
      such that the coordinates of $c_1^\star$ and $c_3^\star$ are $(-1, 0)$ and $(1, 0)$ in $\xi-\eta$ space, respectively
      (see Figure \ref{schematic_2d_p2} (b)).
      The coordinate transformation is given by

     \begin{equation}
     \begin{cases}
     x = \frac{x_3^\star -x_1^\star }{2} \xi  + \frac{ v_3^\star -v_1^\star }{2} \eta + \frac{x_3^\star +x_1^\star }{2}, \\
     v = \frac{ v_3^\star -v_1^\star }{2} \xi  - \frac{ x_3^\star -x_1^\star }{2} \eta + \frac{ v_3^\star + v_1^\star }{2}.
     \end{cases}
     \label{reverse}
     \end{equation}

      By such a transformation, we can get the $\xi-\eta$ coordinate for the point $c_2^\star$ as $(\xi_2, \eta_2)$. Based on the date point $(-1,0)$, $(\xi_2,\eta_2)$ and $(1,0)$, we construct the quadratic curve as follows
         \begin{equation}
         \wideparen{ c_1^\star,c_2^\star, c_3^\star }: \eta(x,v) = \frac{\eta_2}{ \xi_2^2-1 }(  \xi(x,v)^2-1  ).
         \label{parabola}
         \end{equation}

   \end{description}

   \item[(3)] \textbf{\emph{Search algorithm of outer line segments.}}
    We compute all the intersections between grid lines and four curved-sides  of the upstream cell $A_{j}^\star$ and organize them in the  counterclockwise order for outer line segments $\mathcal{L}_q$, $q=1,\cdots,N_o$. Specifically, to find intersection points between grid lines and $\wideparen{ c_1^\star,c_2^\star, c_3^\star }$, we solve the following equations

     \begin{equation}
     \begin{cases}
     x_i =   \frac{x_3^\star -x_1^\star }{2} \xi + \frac{v_3^\star -v_1^\star }{2} \eta
       + \frac{ x_3^\star +x_1^\star }{2} \   (\text{intersection with grid line}\  x=x_i ), \\
      \eta = \frac{\eta_2}{ \xi_2^2-1 }(  \xi^2-1  )
     \end{cases}
     \label{eq:search_x}
     \end{equation}
     and
          \begin{equation}
     \begin{cases}
     v_j  =   \frac{v_3^\star -v_1^\star }{2} \xi - \frac{x_3^\star -x_1^\star }{2} \eta
       + \frac{ v_3^\star +v_1^\star }{2} \   (\text{intersection with grid line}\  v=v_j ) ,\\
      \eta = \frac{\eta_2}{ \xi_2^2-1 }(  \xi^2-1  ) .
     \end{cases}
     \label{eq:search_v}
     \end{equation}
   Using the algorithm provided in Appendix \ref{append:a} to solve \eqref{eq:search_x}-\eqref{eq:search_v},
   we can obtain all intersection points.

   \item[(4)] \textbf{\emph{Evaluation of line integrals along outer line segments.}}
The integral along outer line segments $\sum_{q=1}^{N_o}
\int_{ \mathcal{L}_q } [P  dx +Q  d v ]  $ can be evaluated by the following parameterization for each line segment.
Assume that $\mathcal{L}_q$ is  part of the side $\wideparen{c_1^\star ,c_2^\star, c_3^\star}$.
Substituting \eqref{parabola}  into \eqref{reverse}, we have
     \begin{equation}
     \begin{cases}
     x(\xi) = \frac{x_3^\star -x_1^\star }{2} \xi  + \frac{ v_3^\star -v_1^\star }{2} \frac{\eta_2}{ \xi_2^2-1 }(  \xi^2-1  ) + \frac{x_3^\star +x_1^\star }{2}, \\
     v(\xi) = \frac{ v_3^\star -v_1^\star }{2} \xi  - \frac{ x_3^\star -x_1^\star }{2} \frac{\eta_2}{ \xi_2^2-1 }(  \xi^2-1  ) + \frac{ v_3^\star + v_1^\star }{2}.
     \end{cases}
     \end{equation}

Hence,
\begin{equation}
\label{eq:outerintegral}
\int_{  \mathcal{L}_q } [Pdx+ Qdv ]
= \int_{\xi^{(q)} }^{\xi^{(q+1)} }  [P( x(\xi,\eta),v(\xi,\eta) ) x'(\xi)+ Q( (\xi,\eta),v(\xi,\eta) ) v'(\xi) ]d\xi,
\end{equation}
where $(\xi^{(q)},\eta^{(q)})$ and $(\xi^{(q+1)},\eta^{(q+1)})$  are the start and end points of $\mathcal{L}_q$ in $\xi-\eta$ coordinate, respectively.
The integrand in \eqref{eq:outerintegral} is a polynomial and thus the integral can be exactly computed by a numerical quadrature rule with sufficient degree of precision.

 \end{description}

\subsection{Positivity-preserving limiter}

If the initial condition for the VP system \eqref{vlasov1} is positive, then the solution always stays positive as time evolves.  Such a property is called  positivity preservation. In our SLDG schemes, it can be shown that the updated cell averages at $t^{n+1}$ stay positive, if the numerical solution $f^n$ at $t^n$ is positive.
Similar to \cite{qiu2011positivity,Guo2013discontinuous,guo2015efficient}, in order to preserve positivity of numerical solutions, we further apply a high order positivity-preserving (PP) limiter \cite{zhangshu2010} into the proposed SLDG
scheme, which can be implemented as follows. The numerical solution $f(x,v,t^n)$ in cell $A_j$ is modified by $\widetilde{f}(x,v)$
\begin{equation*}
\widetilde{f}(x,v) = \theta ( f(x,v,t^n)- \overline{f} ) + \overline{f},\
 \theta = \min \left\{  \left| \frac{\overline{f}}{m'-\overline{f}} \right|,1  \right\},
\end{equation*}
where $\overline{f}$ is the cell average of the numerical solution and $m'$ is the minimum value of $f(x,v,t^n)$ over $A_j$.
 For $P^1$ polynomials, the minimum value can be found by comparing the values at four vertices of $A_j$.
 For $P^2$ polynomials, besides the four vertices, all critical points inside $A_j$ should be considered  to determine the function's minimum value. Note that the proposed SLDG schemes with the PP limiter feature the $L^1$ conservation property and hence the $L^1$ stability for nonnegative initial conditions. The proof follows a similar argument in \cite{qiu2011positivity}.

\subsection{The efficient implementation}
\label{section:efficient}

In this subsection, we propose an efficient implementation of the $P^2$ SLDG method with quadratic-curved upstream cells using the third order temporal scheme \eqref{eq:third1}-\eqref{eq:third2}.
As shown in Section \ref{prediction-correction}, for updating the numerical solution from $t^n$ to $t^{n+1}$, this third order scheme includes two prediction steps \eqref{eq: SLDG_t1}, \eqref{eq: SLDG_t2}, and one correction step \eqref{eq: SLDG_t3}. To save some computational cost, we propose to use lower order SLDG schemes in two predictions steps.
The efficient implementation of the $P^2$ SLDG method with quadratic-curved upstream cells using the third order scheme
is summarized in the flow chart named as {\bf Algorithm 1} below. 
Note that the efficient implementation will not compromise that accuracy of the scheme, which can be verified by a simple Taylor expansion. The numerical results presented in the next section also justifies the efficiency of the economical implementation.
%



\bigskip
\fbox{
\begin{minipage}[htb]{0.9\linewidth}
\textbf{Algorithm 1:} The $P^2$ SLDG method with quadratic-curved upstream cells using the third order scheme:
\begin{enumerate}
\item {\bf The first order prediction:}
   \begin{itemize}
     \item Solve the electric field $E$ by the LDG method, based on the solution $f^n$.
     \item Trace the characteristics \eqref{characteristics} for a time step $\Delta t$ by the first order scheme \eqref{eq:first}.
     \item Evolve the solution $f^n$ by using $P^0$ SLDG (i.e. only the cell averages are used and updated) with the quadrilateral approximation to upstream cells to get $f^{n+1,(1)}$.
   \end{itemize}
\item {\bf The second order prediction:}
   \begin{itemize}
     \item Solve the electric field $E$ by the LDG method, based on the solution $f^{n+1,(1)}$.
     \item Trace the characteristics \eqref{characteristics} for a time step $\Delta t$ by the second order scheme \eqref{eq:second}.
     \item Evolve the solution $f^{n+1,(1)}$ by using $P^1$ SLDG with the quadrilateral approximation to upstream cells to get $f^{n+1,(2)}$.
   \end{itemize}
\item {\bf The third order correction:}
   \begin{itemize}
     \item Solve the electric field $E$ by the LDG method, based on the solution $f^{n+1,(2)}$.
     \item Trace the characteristics \eqref{characteristics} for a time step $\Delta t$ by the second order scheme \eqref{eq:third1}-\eqref{eq:third2}.
     \item Evolve the solution $f^{n+1,(2)}$ by using $P^2$ SLDG with quadratic-curved quadrilateral approximation to upstream cells to get $f^{n+1}$.
   \end{itemize}
\end{enumerate}
\end{minipage}
}

\bigskip


%% file: numerical.tex
\section{Numerical Results}
\label{sec3}

In this section, for the VP system, we examine the performance of the proposed SLDG method with second/third order temporal accuracy, denoted by $P^k$ SLDG-time2/3, with quadrilateral or quadratic-curved (QC) quadrilateral approximation to upstream cells (using the notation without or with QC).
We implement the regular as well as efficient versions  of the numerical scheme as described in Section~\ref{section:efficient} (using the notation without or with E). In all of our numerical tests, we let the time step size $\Delta t = \frac{CFL}{ \frac{v_{\max} }{\Delta x} + \frac{\max(|E|)}{\Delta v}}$, where $CFL$ is specified for different runs. For example, $P^2$ SLDG-QC-time3-E-CFL10 refers to the efficient implementation of the SLDG scheme with $P^2$ solution space, with quadratic-curved quadrilateral approximation to upstream cells, using third order scheme in characteristics tracing and with $CFL=10$. We apply PP limiter for all test examples. We also note that the proposed SLDG methods are stable and highly accurate under very large $CFL$ numbers as shown in numerical tests in this section. However, if an excessively large $CFL$ number is used,  some approximated upstream cells may become distorted, leading to a breakdown of the scheme. When the distortion happens in the simulation, a smaller time step should be used so that no distortion appears.

In this section, we demonstrate the following different aspects via extensive numerical tests of the proposed algorithm on a set of benchmark VP examples.
\begin{enumerate}
\item {\em Using a high order characteristics tracing scheme.} For weak Landau damping, we benchmark the numerical damping rate of the electrostatic field $E$  against the theoretical value from the linear theory. In particular, we test the SLDG method with second and third order characteristics tracing schemes using $CFL$ numbers as large as $20$. In Figure~\ref{weakLandau_e2}, we showed the advantage of using a third order characteristics tracing scheme (compared with a second order one) for its superior performance in capturing the correct damping rate with $CFL$ as large as $20$.
\item {\em Spatial order of convergence: the need to use quadratic curves in approximating sides of upstream cells.} We test the spatial order of convergence in Table~\ref{spatial_1} for strong Landau damping, and in Table~\ref{spatial_two1} for two stream instabilities.  The computational effort in using a quadratic-curved quadrilateral approximation of upstream cells is justified by smaller error magnitudes observed and the corresponding third order convergence for the SLDG scheme with $P^2$ polynomial space.
\item {\em Temporal order of convergence.} We test the temporal order of convergence in Table~\ref{temporal_2} for strong Landau damping and in Table~\ref{temporal_two1_reverse} and \ref{temporal_two1} for two stream instabilities. Simulations with large $CFL$s, ranging from $5$ to $50$ or more, provide decent results.
\item {\em CPU savings for the efficient implementation of the SLDG scheme.} We compare the numerical performance and CPU time for the ``regular" and ``efficient" implementation of the SLDG scheme in Table~\ref{spatial_1} and \ref{spatial_two1}. In particular, for the $P^2$ SLDG scheme, more than $30\%$ savings in CPU time are observed, while the accuracy is not compromised. Due to similar performance and significant savings in CPU time, most of our tests are done base on the efficient implementation.
\item {\em CPU comparison between the SLDG scheme and the classical RKDG.} We perform CPU comparison between the SLDG and RKDG methods in Tables \ref{temporal_2} and \ref{temporal_two1_reverse}. Per time step evolution, an SLDG scheme takes about four times as much CPU time as that of a RKDG counterpart with the same order accuracy. On the other hand, for stability of an Eulerian RKDG method, the upper bound of the $CFL$ is about $\frac{1}{2k+1}$, with $k$ being the polynomial degree; while for an SLDG scheme, the $CFL$ number can be taken as large as $50$, leading to tremendous savings in CPU time ($50\%$ to $90\%$ savings). We can choose the time stepping size purely for accuracy consideration, without much constraint from stability.
\item {\em Preservation of mass and other physical norms.}
In the VP system, the following physical quantities should remain constant in time. Tracking relative deviations of these quantities numerically provides a good measurement of the quality of numerical schemes.
Our proposed SLDG scheme is locally and globally mass conservative.
We will show comparable (sometimes superior) performance  in preserving the physical norms for the proposed SLDG scheme with large $CFL$s.
\begin{description}
  \item[1.] Mass:
  \begin{equation*}
  \text{Mass} = \int_v\int_x f(x,v,t)dxdv.
  \end{equation*}
  \item[2.] $L^p$ norm, $1\leq p<\infty$:
  \begin{equation*}
  \| f \|_p = \left( \int_v\int_x |f(x,v,t)|^p dxdv\right)^{\frac{1}{p} }.
  \end{equation*}
  \item[3.] Energy:
  \begin{equation*}
   \text{Energy} = \int_v\int_x f(x,v,t) v^2 dxdv + \int_x E^2(x,t) dx,
  \end{equation*}
  where $E(x,t)$ is the electric field.
  \item[4.] Entropy:
  \begin{equation*}
  \text{Entropy} = \int_v\int_x f(x,v,t) \log( f(x,v,t) ) dxdv.
  \end{equation*}
\end{description}
\end{enumerate}


\begin{exa}
\emph{(Weak Landau damping.)}
Consider weak Landau damping for the VP system. The initial condition is set to be the following perturbed equilibrium
\begin{equation}
f(x,v,t=0) = \frac{ 1 }{\sqrt{2\pi} } ( 1+ \alpha \cos(kx)  ) \exp \left( -\frac{v^2}{2} \right) ,
\label{weaklandau_init}
\end{equation}
with $\alpha = 0.01$ and $k=0.5$. Our computational domain is $[0,4\pi]\times[-v_{\max},v_{\max}]$.
We truncate the velocity domain at $v_{\max}=2\pi$.
This test case has been numerically investigated by several authors (e.g. see \cite{FilbetSB, Qiu_Christlieb,rossmanith2011positivity,heath2012discontinuous,guo2013hybrid,gucclu2014arbitrarily,cai2016conservative,qiu2017high}).

In Figure \ref{weakLandau_e2}, we present the time evolution of $L^2$ norm of the electric field (in semi-log scale) for $P^2$ SLDG-time2-E and $P^2$ SLDG-time3-E schemes using a mesh of $128\times128$ elements and different $CFL$s.
As expected, the decay rate from simulations with $CFL=1$ matches well with the theoretical value $\gamma=-0.1533$ \cite{FilbetSB} (the solid line in the same plots). For $CFL=10$, both results match well with the theoretical value; for $CFL=20$, the third order characteristics tracing scheme exhibits superior performance, compared with the second order one, in capturing the correct damping rate in the long run.


\begin{figure}[h!]
\centering
\includegraphics[height=60mm]{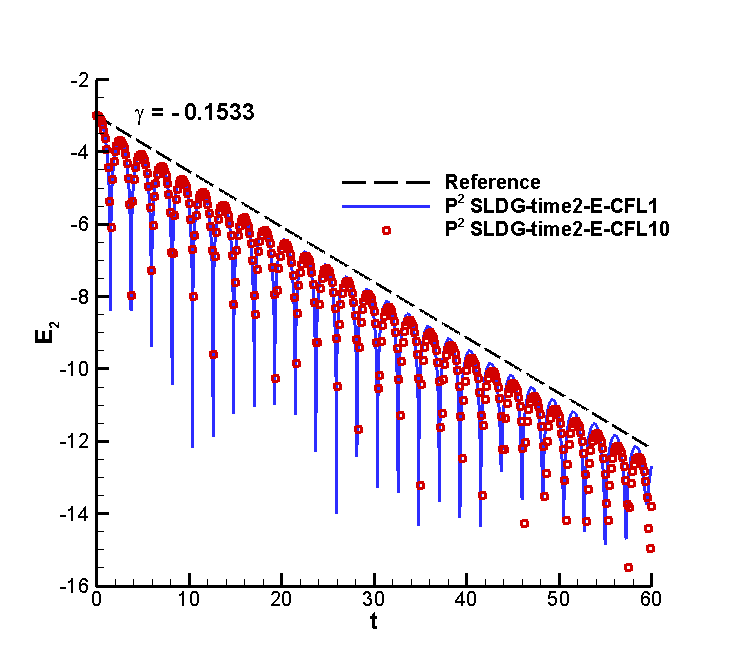}
\includegraphics[height=60mm]{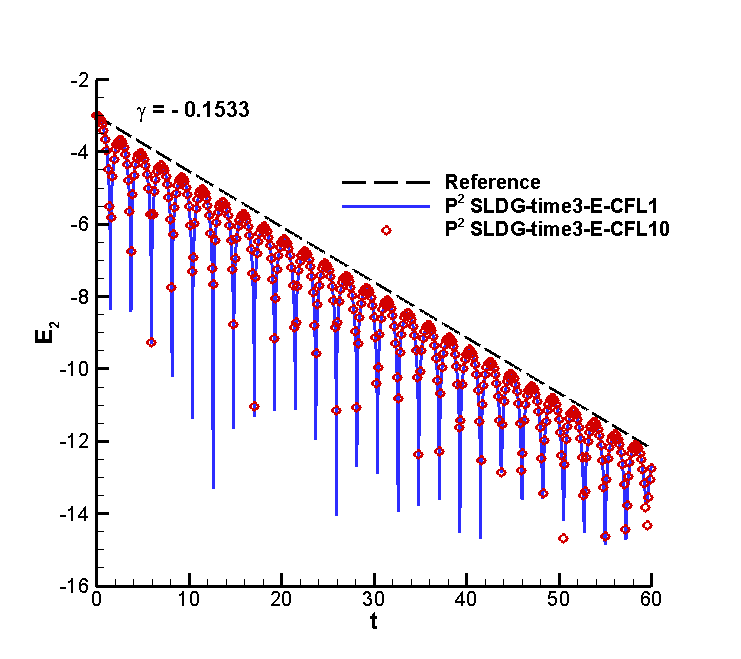}
\includegraphics[height=60mm]{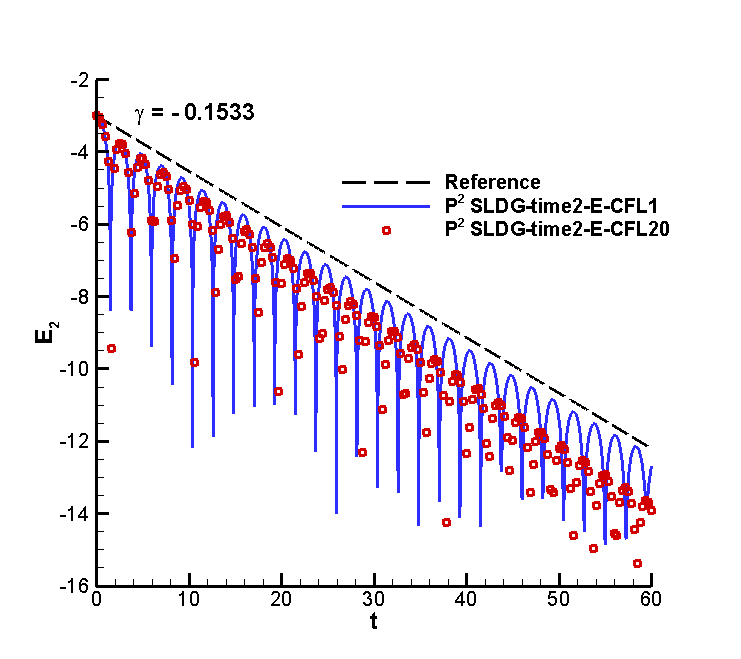}
\includegraphics[height=60mm]{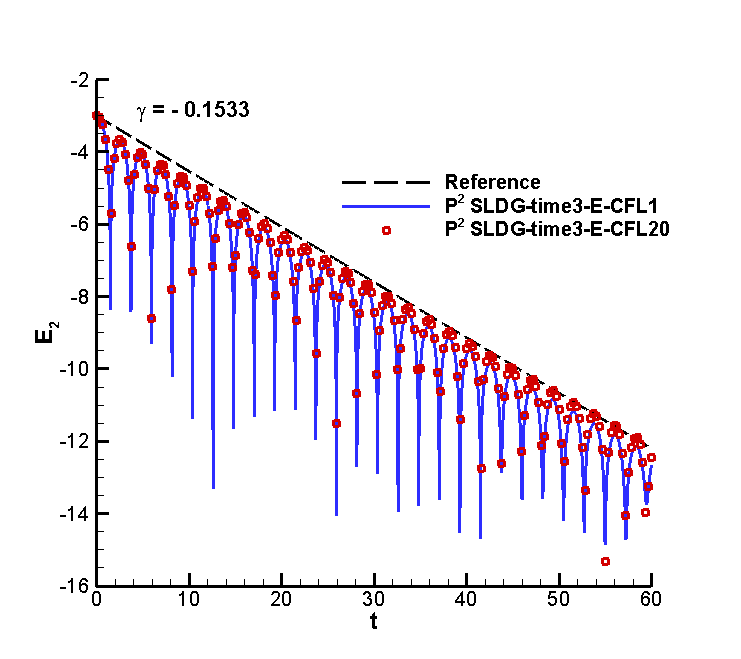}
\caption{Weak Landau damping. Time evolution of electric field in $L^2$ norm.
Solid line: $CFL=1$. Square: $CFL=10$ or $CFL=20$.
The 'Reference' line (dashed line) corresponds to the exponential decay in the amplitude of the oscillation ($\gamma=-0.1533$).}
\label{weakLandau_e2}
\end{figure}
\end{exa}

\begin{exa}
\emph{(Strong Landau damping.)} Consider strong Landau damping for the VP system.  The initial condition is the same as weak one \eqref{weaklandau_init},
but with a larger perturbation parameter $\alpha = 0.5$. The computational domain is $[0,4\pi]\times[-2\pi,2\pi]$.

In Table \ref{spatial_1}, we test the spatial convergence of the proposed SLDG scheme with the third order characteristic tracing scheme. We set $CFL=0.1$ to minimize the error from time discretization. The well-known time reversibility of the VP system is used to test the order of convergence.
In particular, one can integrate the VP system forward to some time $T$, and then reverse
the velocity field of the solution and continue to integrate the system by the same amount of time $T$. Then, the solution should recover the initial condition with reverse velocity field, which can be used as a reference solution. We show the $L^2$ and $L^\infty$ errors and the corresponding orders of convergence for $P^k$ SLDG-(QC)-time3-(E) schemes with $CFL=0.1$ in Table \ref{spatial_1}. Second order convergence is observed for $P^1$ SLDG scheme as expected. Second order convergence, with smaller error magnitudes than those for $P^1$ SLDG scheme, is observed for $P^2$ SLDG scheme with quadrilateral approximation to upstream cells. Such second order convergence is due to the use of straight lines in approximating sides of upstream cells, even though $P^2$ approximation space is employed. More importantly, if quadratic curves are used in approximating sides of upstream cells, the error is further reduced and the order of convergence is improved to third order, see the results for $P^2$ SLDG-QC scheme.  $30\%-40\%$ savings in CPU time are observed  for the {\em efficient} implementation of $P^2$ SLDG schemes when compare with the CPU time needed for the {\em regular} implementation, while the accuracy is not compromised.


Second, temporal convergence from the characteristics tracing scheme used in the proposed SLDG method is tested.
Table \ref{temporal_2} summarizes the $L^2$ and $L^\infty$ errors and the corresponding temporal convergence rates for the $P^k$ SLDG methods with the second and third order characteristic tracing schemes, where the time reversibility property of the VP system is used to compute the error. To make the temporal error dominant, we use a spatial mesh of $160\times160$ elements.
In Table \ref{temporal_2}, around second and third order temporal accuracy is observed, for $P^1$ SLDG-time2-E and $P^2$ SLDG-QC-time3-E schemes, respectively, with $CFL$s ranging from $5$ to $25$. Also in Table \ref{temporal_2}, we perform CPU time comparison between SLDG and RKDG schemes. Per time step, the SLDG scheme costs about $4$ times as much CPU time as that of the classical RKDG method. On the other hand, the SLDG scheme allows for a much larger $CFL$. For example,  if $CFL=25$, then time stepping size of the SLDG scheme is $125$ times as large as that for a RKDG scheme, leading to tremendous savings in CPU time. Note that, the RKDG code, that we use for CPU comparison, has been optimized for its CPU efficiency.

We show the time evolution of the electric field in the $L^2$ norm (in semi-log scale) in Figure \ref{strongLandau_e2}. The linear decay rate $\gamma_1$ (measured as the slope of a line originating from the local maximum of the second peak to the third peak), as well as the growth rate $\gamma_2$ (measured as the slope of a line originating from the local maximum of the tenth peak to the sixteenth peak) are summarized in Table \ref{Landau_gamma}, and they agree with the results reported in the literature
\cite{FilbetSB, rossmanith2011positivity, heath2012discontinuous, guo2013hybrid}. Furthermore, in Figure \ref{strongLandau}, we plot time evolutions of the relative derivation of the discrete $L^1$ norm, $L^2$ norm, energy and entropy.
In particular, we observe that (1)
 The error for the $L^1$ norm (on the order of $10^{-9}$) is due to the truncation of the velocity domain, which can be further reduced by using a larger velocity domain in simulations; (2) In general,
 the $P^2$ SLDG method does a better job in conserving these physical norms than the  $P^1$ SLDG method; (3) Compared to the SLDG schemes with larger $CFL$s, the schemes with smaller $CFL$s are able to better conserve the energy, but perform worse in conserving  the $L^2$ norm and entropy.
 (4) The SLDG methods outperform the RKDG method in conserving energy.
 In Figure \ref{strongLandau1}, we present the contour plots of the solutions at $T=40$ computed by the  $P^2$ SLDG-QC-time3-E scheme method and the $P^2$ RKDG method with the mesh of $160\times160$ elements.
We observe that the SLDG scheme with $CFL=10$ and the RKDG scheme with $CFL=0.2$ generate very consistent numerical results. Meanwhile, the SLDG method with $CFL=30$ is still stable and able to generate decent results: the main structures of the solution are captured,  but some mild wiggles are observed.

 \begin{table}[!ht]\small
\caption{Strong Landau damping with $T=0.5$. Order of accuracy in space and CPU time for $P^k$ SLDG-(QC)-time3-(E) scheme. We set $CFL=0.1$ so that the spatial error is the dominant error.}
\vspace{0.1in}
\centering
\begin{tabular}{c cc  cc c  cc cc c}
\hline

{ Mesh}   &{$L^2$ error} & Order    &{$L^\infty$ error} & Order &CPU
 &{$L^2$ error} & Order    &{$L^\infty$ error} & Order &CPU\\
\hline
& \multicolumn{5}{c}{$P^1$ SLDG-time3}   & \multicolumn{5}{c}{$P^1$ SLDG-time3-E} \\
 \cmidrule(lr){2-6} \cmidrule(lr){7-11}
  $32^2$  &  1.21E-3 &        &   1.18E-2 &         & 2.28 &   1.21E-3 & &     1.18E-2 &  & 2.17\\
  $64^2$  &  3.15E-4 &     1.94 &     3.49E-3 &     1.76 & 17.51 &   3.15E-4 &     1.94 &     3.49E-3 &     1.76 & 16.09 \\
  $96^2$  &  1.42E-4 &     1.97 &     1.61E-3 &     1.91 & 58.50 &   1.42E-4 &     1.97 &     1.61E-3 &     1.91 & 55.10 \\
  $128^2$ &  8.02E-5 &     1.98 &     9.18E-4 &     1.95 & 139.60 &   8.02E-5 &     1.98 &     9.18E-4 &     1.95 & 129.78\\
  $160^2$ &  5.15E-5 &     1.99 &     5.92E-4 &     1.97 & 275.96 &   5.15E-5 &     1.99 &     5.92E-4 &     1.97 & 252.92\\
\hline
& \multicolumn{5}{c}{$P^2$ SLDG-time3}   & \multicolumn{5}{c}{$P^2$ SLDG-time3-E} \\
 \cmidrule(lr){2-6} \cmidrule(lr){7-11}
  $32^2$  &  2.18E-4 & &     1.97E-3 &         & 4.96 &   2.18E-4 & &    1.95E-3 &  & 3.12\\
  $64^2$  &  5.57E-5 &     1.97 &     5.03E-4 &     1.97 & 37.50 &   5.57E-5 &    1.97 &    4.72E-4 &     2.04 & 24.68 \\
  $96^2$  &  2.50E-5 &     1.98 &     2.25E-4 &     1.99 & 126.96 &   2.50E-5 &     1.98 &     2.11E-4 &     1.99 & 81.09 \\
  $128^2$ &  1.41E-5 &     1.98 &     1.27E-4 &     2.00 & 315.18 &  1.41E-5 &     1.98 &     1.19E-4 &     2.00& 201.04\\
  $160^2$ &  9.07E-6 &     1.99 &     8.10E-5 &     2.00 & 602.12 &  9.07E-6 &     1.99 &     7.60E-5 &     2.00& 391.17\\
\hline
& \multicolumn{5}{c}{$P^2$ SLDG-QC-time3}   & \multicolumn{5}{c}{$P^2$ SLDG-QC-time3-E} \\
 \cmidrule(lr){2-6} \cmidrule(lr){7-11}
  $32^2$  &  8.32E-5 & &     1.08E-3 &         & 5.29 &   8.32E-5 & &     1.08E-3 &  & 3.59\\
  $64^2$  &  1.02E-5 &     3.03 &     1.38E-4 &     2.97 & 41.50 &   1.02E-5 &     3.03 &     1.36E-4 &     2.98 & 27.60 \\
  $96^2$  &  3.00E-6 &     3.02 &     4.08E-5 &     3.00 & 141.18 &   3.00E-6 &     3.02 &     4.02E-5 &     3.01 & 91.87 \\
  $128^2$ &  1.26E-6 &     3.01 &     1.71E-5 &     3.03 & 334.73 &   1.26E-6 &     3.01 &     1.68E-5 &     3.02 & 221.90\\
  $160^2$ &  6.43E-7 &     3.01 &     8.58E-6 &     3.08 & 645.90 &   6.43E-7 &     3.01 &     8.45E-6 &     3.09 & 433.79\\
\hline
\end{tabular}
\label{spatial_1}
\end{table}

 \begin{table}[!ht]\small
\caption{Strong Landau damping with $T=0.5$, and with the mesh of $160\times160$. Temporal order of convergence and CPU comparison between SLDG and RKDG schemes.  }
\vspace{0.1in}
\centering
\begin{tabular}{c cc  cc  c }
\hline

{ $CFL$}  &{$L^2$ error} & Order    &{$L^\infty$ error} & Order & CPU (sec)   \\
\hline
 \multicolumn{6}{l}{$P^1$ SLDG-time2-E}  \\
 0.3 &        5.05E-05 & --&     5.90E-04 &    -- & {\textcolor{red}{47.79}}\\
   5 &         3.97E-05 & -- &     5.67E-04 &      --         &{\textcolor{red}{3.23}}\\
   10 &        4.67E-05 &     0.23 &     6.27E-04 &     0.15 &{\textcolor{red}{1.67}}\\
   15 &        1.04E-04 &     1.97 &     1.04E-03 &     1.24  &{\textcolor{red}{1.32}}\\
   20 &       1.96E-04 &     2.21 &     1.44E-03 &     1.14 &{\textcolor{red}{1.01}}\\
   25 &        3.91E-04 &     3.08 &     2.45E-03 &     2.37 &{\textcolor{red}{0.71}}\\
   \hline
    \multicolumn{6}{l}{$P^1$ RKDG}  \\
   0.3 &         4.36E-05 &    -- &     4.17E-04 &    -- &  {\textcolor{red}{5.26}}\\
   \hline
\multicolumn{6}{l}{$P^2$ SLDG-QC-time3-E }  \\
  0.2 &        6.37E-07 & -- &     8.43E-06 &  --  & {\textcolor{red}{229.92}}\\
   5 &         2.35E-06 & -- &     1.30E-05 &   --                &{\textcolor{red}{9.50}}\\
   10 &        1.69E-05 &     2.85 &     9.87E-05 &     2.92 &{\textcolor{red}{5.14}}\\
   15 &        6.15E-05 &     3.19 &     2.83E-04 &     2.60 &{\textcolor{red}{3.89}}\\
   20 &        1.32E-04 &     2.65 &     5.97E-04 &     2.60 &{\textcolor{red}{2.98}}\\
   25 &        2.43E-04 &     2.74 &     1.11E-03 &     2.79 &{\textcolor{red}{2.01}}\\
   \hline
   \multicolumn{6}{l}{$P^2$ RKDG}  \\
   0.2 &        8.68E-07 &    -- &     1.09E-05 &    -- & {\textcolor{red}{41.39}}\\
\hline
\end{tabular}
\label{temporal_2}
\end{table}



\begin{table}[!ht]\small
\caption{Strong Landau damping. $T=0.5$. A mesh of $160\times160$ is used. The linear decay rate $\gamma_1$ and the growth rate $\gamma_2$. }
\vspace{0.1in}
\centering
\begin{tabular}{c cc  cc }
\hline
{ $CFL$}  &$\gamma_1$& $\gamma_2$   &$\gamma_1$& $\gamma_2$   \\
\hline
 &\multicolumn{2}{c}{$P^1$ SLDG-time2-E}  &\multicolumn{2}{c}{$P^2$ SLDG-QC-time3-E}  \\
 \cmidrule(lr){2-3} \cmidrule(lr){4-5}
 1 & -0.2907 & 0.0847   &  -0.2907 & 0.0868 \\
 10 & -0.2875  & 0.0848 & -0.2842 & 0.0866  \\
 20 & -0.3054  & 0.0840 & -0.3125 & 0.0867 \\
\hline
\end{tabular}
\label{Landau_gamma}
\end{table}

\begin{figure}[h!]
\centering
\includegraphics[height=65mm]{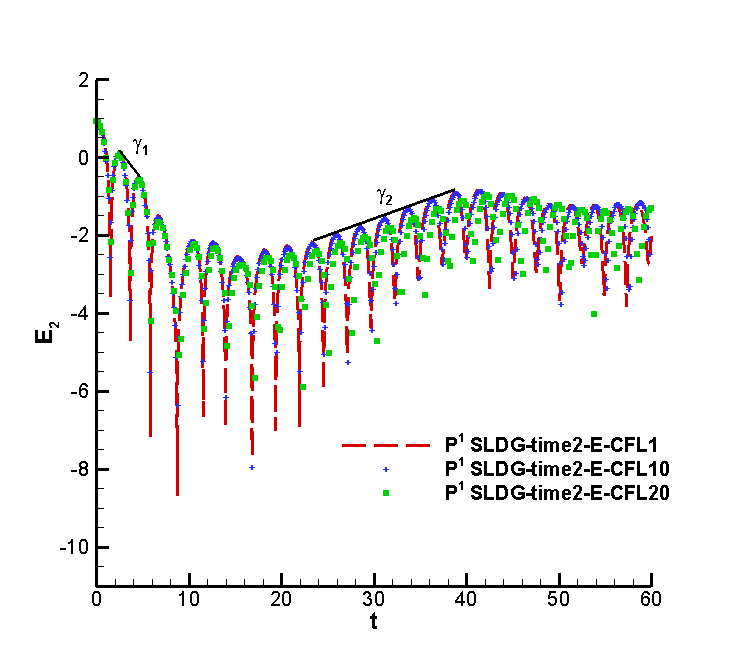}
\includegraphics[height=65mm]{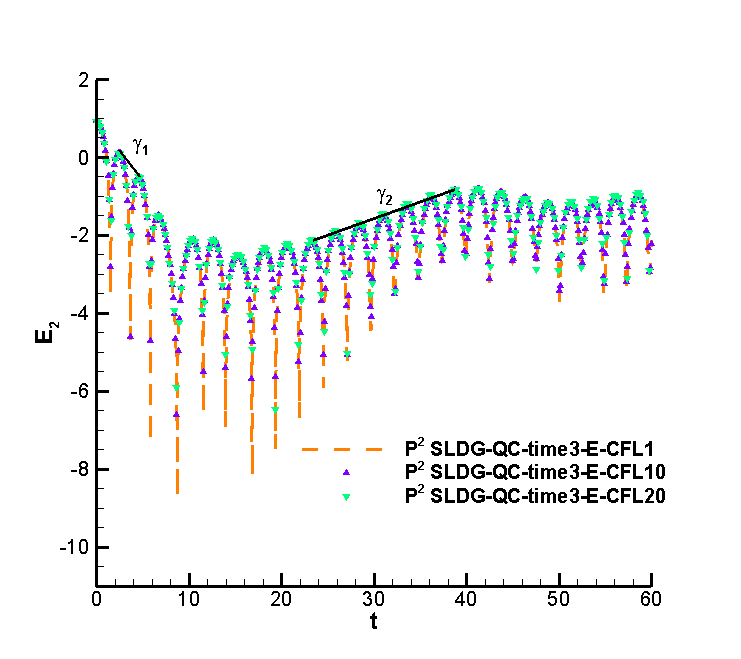}
\caption{Strong Landau damping: The SLDG schemes are equipped with the PP limiter. Time evolution of the electric field in $L^2$. }
\label{strongLandau_e2}
\end{figure}

\begin{figure}[h!]
\centering
\includegraphics[height=65mm]{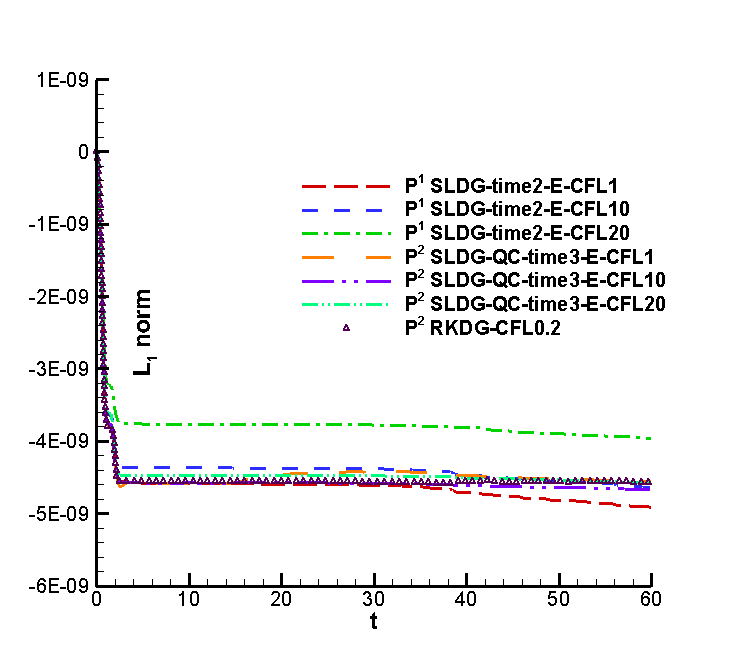}
\includegraphics[height=65mm]{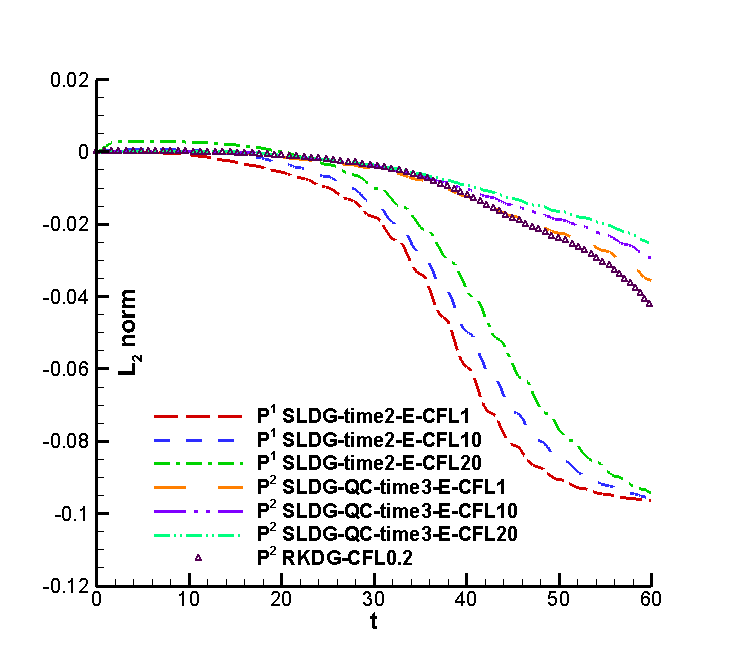}
\includegraphics[height=65mm]{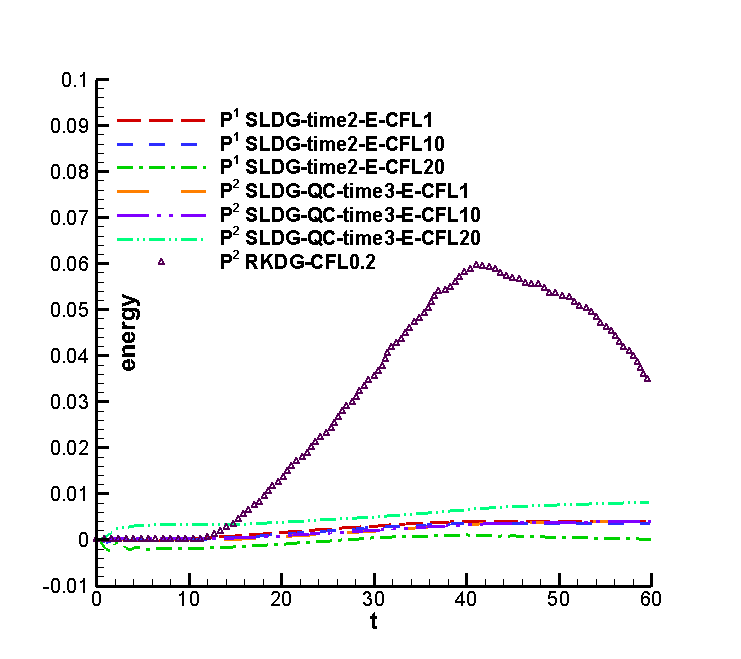}
\includegraphics[height=65mm]{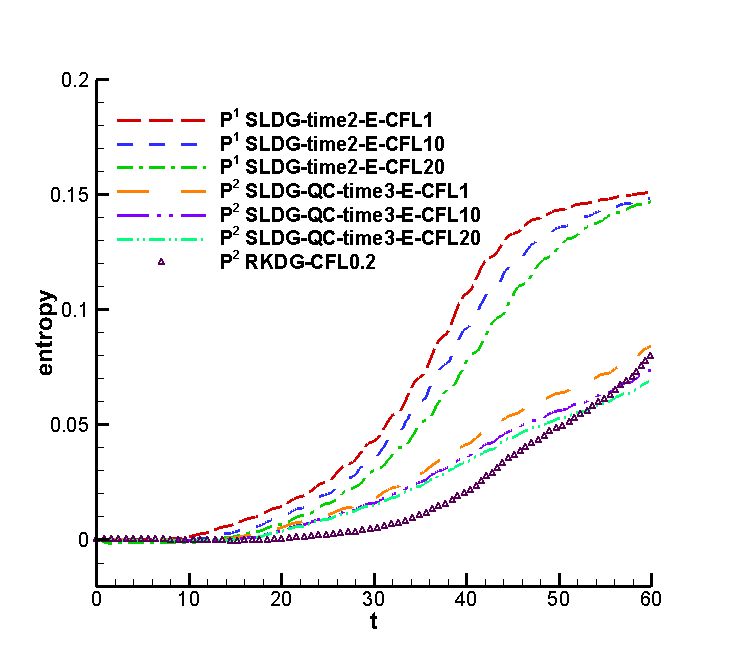}
\caption{Strong Landau damping. Time evolution of  the relative deviations of $L^1$ (upper left) and $L^2$ (upper right) norms of the solution as well as the discrete kinetic energy (lower left) and entropy (lower right).}
\label{strongLandau}
\end{figure}

\begin{figure}[h!]
\centering
\subfigure[$P^2$ SLDG-QC-time3-E with $CFL=10$]{
\includegraphics[height=45mm]{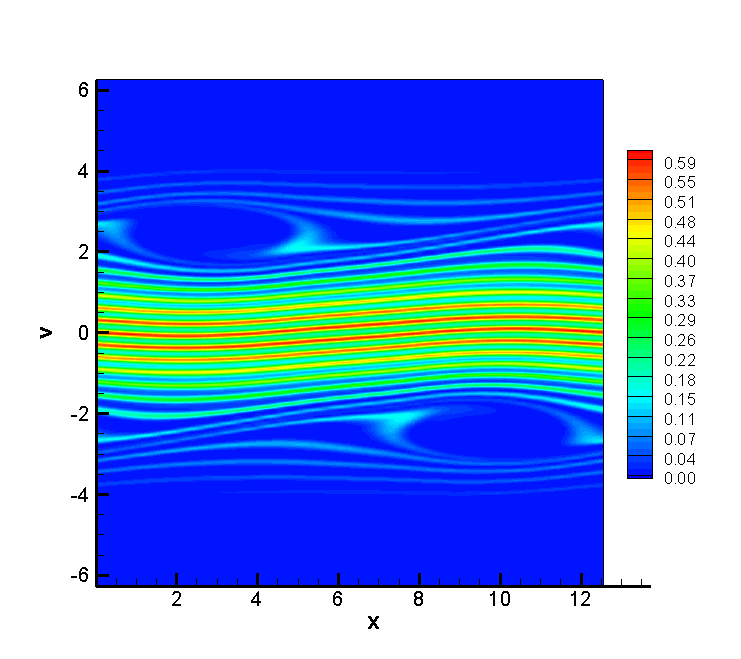} }
\subfigure[$P^2$ SLDG-QC-time3-E with $CFL=30$]{
\includegraphics[height=45mm]{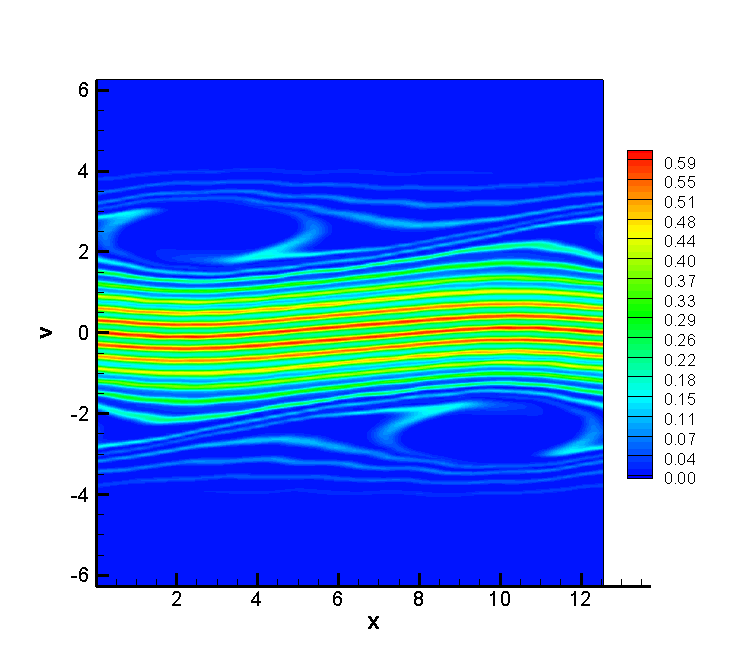} }
\subfigure[$P^2$ RKDG with $CFL=0.2$]{
\includegraphics[height=45mm]{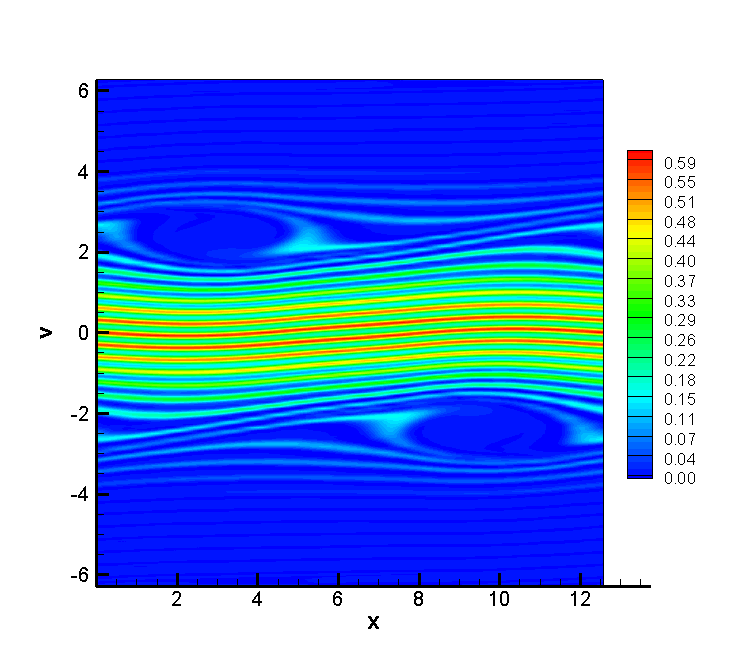} }
\caption{Strong Landau damping with the spatial mesh of $160\times160$. $T=40$.
(a): $P^2$ SLDG-QC-time3-E with $CFL=10$. (b): $P^2$ SLDG-QC-time3-E with $CFL=30$.
(c): $P^2$ RKDG with $CFL=0.2$.}
\label{strongLandau1}
\end{figure}
\end{exa}

\begin{exa}
\label{exa:two1}
\emph{(Two Stream instability I.)}
Consider two stream instability, with an unstable initial distribution function,
\begin{equation*}
f(x,v,t=0) = \frac{2}{7\sqrt{2\pi} } (1+5v^2) ( 1 + \alpha \left(  ( \cos(2kx)+\cos(3kx) )/1.2 + \cos(kx) \right) )\exp\left(-\frac{v^2}{2} \right)
\end{equation*}
with $\alpha =0.01,\ k=0.5$ on the computation domain $[0,4\pi]\times[-10,10]$.

We first present the spatial convergence of the SLDG schemes in Table \ref{spatial_two1}. As in the previous example, we make use of the time reversibility of VP system to test the order of convergence. Slightly less than second order accuracy is observed for the $P^1$ SLDG-time2-E scheme.
The error magnitude of the $P^2$ SLDG-time3-E scheme is much smaller than the corresponding one of the $P^1$  SLDG-time2-E scheme, while order reduction is still observed due to the fact that we use straight lines to approximate sides of upstream cells. When quadratic curves are used to approximate sides of upstream cells, the errors are further reduced and around third order convergence is observed for the $P^2$ SLDG-QC scheme. Similar to strong Landau damping, $30\%-40\%$ savings in CPU time are observed for $P^2$ SLDG-(QC)-time3 schemes with {\em efficient} implementation.

We test the temporal convergence of the SLDG schemes for this example. In order to make the temporal error dominant, we use a spatial mesh of $160\times160$ elements. The convergence results, as well as CPU comparison with RKDG methods, are presented in Table \ref{temporal_two1_reverse}.
For this example, even from the accuracy point of view,  $CFL$ for the $P^1$ SLDG scheme can be taken to be as large as $50$ (compared with $0.3$ for a $P^1$ RKDG scheme) with comparable performance in terms of magnitude of errors. Again, the corresponding savings in CPU time are tremendous compared with RKDG methods: it takes $7.35$ seconds for the RKDG method to reach the final integration time $T=0.5$; while it takes less than one second for the SLDG method. Similar comments apply to the comparison of CPU time between the $P^2$ SLDG and RKDG schemes. Table \ref{temporal_two1_reverse} demonstrates the superior performance and efficiency of the proposed SLDG method. In order to show the temporal convergence rate of the SLDG scheme, we present Table \ref{temporal_two1}, where the reference solution comes from the simulation with the same spatial mesh but with relatively small $CFL=0.1$. Expected second and third order convergence rates are observed.



Lastly, we show time evolution of the electric field in the $L^2$ norm (in semi-log scale) for the $P^1$ SLDG-time2-E and $P^2$ SLDG-QC-time3-E schemes using a mesh of $160\times160$ elements and different $CFL$s in Figure \ref{two1_e2}.
In Figure \ref{two1norm},
we show the relative derivation of the discrete $L^1$ norm, $L^2$ norm, energy and entropy.
We observe that
all methods are able to conserve the $L^1$ norm up to the truncation error from the velocity domain. The ability of SLDG methods to conserve these physical norms is satisfactory and comparable to the RKDG method.
In Figure \ref{two1solution}, we plot the numerical solutions of phase space profiles  at $T=40$.
We observe that the solutions computed by the  $P^2$  SLDG-QC-time3-E scheme with $CFL=10,30$ are consistent with that by the $P^2$ RKDG scheme with $CFL=0.2$.


 \begin{table}[!ht]\small
\caption{Two Stream instability I at $T=0.5$. Spatial order of accuracy and CPU time for the SLDG method. We set $CFL=0.1$ so that the spatial error is the dominant error.}
\vspace{0.1in}
\centering
\begin{tabular}{c cc  cc c  cc cc c}
\hline

{ Mesh}   &{$L^2$ error} & Order    &{$L^\infty$ error} & Order &CPU
 &{$L^2$ error} & Order    &{$L^\infty$ error} & Order &CPU\\
\hline
& \multicolumn{5}{c}{$P^1$ SLDG-time3}   & \multicolumn{5}{c}{$P^1$ SLDG-time3-E} \\
 \cmidrule(lr){2-6} \cmidrule(lr){7-11}
  $32^2$  &  4.28E-3 & &     2.28E-2 &         & 3.17 &   4.28E-3 & &     2.28E-2 &  & 2.93\\
  $64^2$  &  1.10E-3 &     1.97 &     7.29E-3 &     1.65 & 24.89 &   1.10E-3 &     1.97 &     7.29E-3 &     1.65 & 23.04\\
  $96^2$  &  4.91E-4 &     1.98 &     3.60E-3 &     1.74 & 83.18 &   4.91E-4 &     1.98 &     3.60E-3 &     1.74 & 75.06\\
  $128^2$ &  2.79E-4 &     1.96 &     2.21E-3 &     1.69 & 194.60&  2.79E-4 &     1.96 &     2.21E-3 &     1.69 &182.12\\
  $160^2$ &  1.81E-4 &     1.93 &     1.53E-3 &     1.65 & 395.15&  1.81E-4 &     1.93 &     1.53E-3 &     1.65 &354.96\\
\hline
& \multicolumn{5}{c}{$P^2$ SLDG-time3}   & \multicolumn{5}{c}{$P^2$ SLDG-time3-E} \\
 \cmidrule(lr){2-6} \cmidrule(lr){7-11}
  $32^2$  &  5.33E-4 & &     2.80E-3 &         & 6.78 &   5.33E-4 & &     2.79E-3 &  & 4.35\\
  $64^2$  &  7.47E-5 &     2.84 &     3.88E-4 &     2.85 & 53.62 &   7.47E-5 &     2.84 &     3.87E-4 &     2.85 & 34.25\\
  $96^2$  &  2.52E-5 &     2.67 &     1.63E-4 &     2.14 & 179.95 &   2.52E-5 &    2.68 &     1.63E-4 &     2.14 & 113.35\\
  $128^2$ &  1.30E-5 &     2.31 &     1.20E-4 &     1.08 & 419.73 &  1.30E-5 &     2.31 &     1.19E-4 &     1.08& 270.32\\
  $160^2$ &  8.61E-6 &     1.84 &     1.00E-4 &     0.80 & 821.60 &  8.60E-6 &     1.84 &     9.99E-5 &     0.80& 520.23\\
\hline
& \multicolumn{5}{c}{$P^2$ SLDG-QC-time3}   & \multicolumn{5}{c}{$P^2$ SLDG-QC-time3-E} \\
 \cmidrule(lr){2-6} \cmidrule(lr){7-11}
  $32^2$  &  5.34E-4 & &     2.81E-3 &         & 7.45 &   5.34E-4 & &     2.81E-3 &  & 4.92\\
  $64^2$  &  7.41E-5 &     2.85 &     3.72E-4 &     2.92 & 58.42 &   7.41E-5 &     2.85 &     3.72E-4 &     2.92 & 38.93 \\
  $96^2$  &  2.39E-5 &     2.79 &     1.28E-4 &     2.63 & 202.46 &   2.39E-5 &     2.79 &     1.28E-4 &   2.63 & 131.04 \\
  $128^2$ &  1.08E-5 &     2.76 &     6.03E-5 &     2.62 & 469.59 &   1.08E-5 &     2.76 &     6.03E-5 &    2.62 & 312.01\\
  $160^2$ &  5.87E-6 &     2.74 &     3.29E-5 &     2.72 & 891.53 &   5.87E-6 &     2.74 &     3.29E-5 &     2.72 & 596.71\\
\hline
\end{tabular}
\label{spatial_two1}
\end{table}

\begin{table}[!ht]\small
\caption{Two Stream instability I at $T=0.5$. A mesh of $160\times160$ is used. Temporal order of convergence for the SLDG method via the time
reversibility of the VP system.  }
\vspace{0.1in}
\centering
\begin{tabular}{c cc  cc  c }
\hline

{ $CFL$}  &{$L^2$ error} & Order    &{$L^\infty$ error} & Order & CPU   \\
\hline
 \multicolumn{6}{l}{$P^1$ SLDG-time2-E}  \\
 0.3 &        1.81E-04 &-- &     1.53E-03 &    -- & {\textcolor{red}{64.65}}\\
%
   45 &       2.07E-04 &     -- &     1.56E-03 &     --  &{\textcolor{red}{0.70}}\\
   50 &       2.31E-04 &     1.05 &     1.63E-03 &     0.42 &{\textcolor{red}{0.70}}\\
   55 &       2.73E-04 &     1.76 &     1.91E-03 &     1.65 &{\textcolor{red}{0.70}}\\
   60 &       3.42E-04 &     2.60 &     2.36E-03 &     2.46  &{\textcolor{red}{0.70}}\\
   65 &       4.21E-04 &     2.58 &     2.87E-03 &     2.42 &{\textcolor{red}{0.39}}\\
   \hline
    \multicolumn{6}{l}{$P^1$ RKDG}  \\
   0.3 &        1.83E-04 &    -- &     1.51E-03 &     -- & {\textcolor{red}{7.35}}\\
   \hline
\multicolumn{6}{l}{$P^2$ SLDG-QC-time3-E }  \\
  0.2 &        5.87E-06 & -- &     3.29E-05 &  --  & {\textcolor{red}{307.57}}\\
   5 &         5.85E-06 & -- &     3.25E-05 &   --         & {\textcolor{red}{12.31}}\\
   10 &        5.96E-06 &     0.03 &     3.30E-05 &     0.02 &{\textcolor{red}{6.75}}\\
   15 &        7.23E-06 &     0.48 &     4.31E-05 &     0.66 &{\textcolor{red}{4.81}}\\
   20 &        1.18E-05 &     1.69 &     7.70E-05 &     2.01 &{\textcolor{red}{3.93}}\\
   25 &        1.89E-05 &     2.13 &     1.23E-04 &     2.11 &{\textcolor{red}{3.00}}\\

   30 &        3.52E-05 &     3.40 &     2.32E-04 &     3.46& {\textcolor{red}{2.98}}\\
   35 &        4.72E-05 &     1.91 &     3.17E-04 &     2.03 &{\textcolor{red}{2.01}}\\
   40 &        6.49E-05 &     2.38 &     4.23E-04 &     2.16 &{\textcolor{red}{2.00}}\\
   45 &        9.75E-05 &     3.45 &     6.16E-04 &     3.19 &{\textcolor{red}{1.97}}\\
   50 &        1.46E-04 &     3.82 &     9.28E-04 &     3.90 &{\textcolor{red}{1.92}}\\
   \hline
   \multicolumn{6}{l}{$P^2$ RKDG}  \\
   0.2 &        6.98E-06 &   --  &     3.43E-05 &    -- & {\textcolor{red}{57.29}}\\
\hline
\end{tabular}
\label{temporal_two1_reverse}
\end{table}

 \begin{table}[!ht]\small
\caption{Two Stream instability I at $T=0.5$. A mesh of $160\times160$ is used. Temporal order of convergence for the SLDG method.  }
\vspace{0.1in}
\centering
\begin{tabular}{c cc  cc  }
\hline

{ $CFL$}  &{$L^2$ error} & Order    &{$L^\infty$ error} & Order     \\
\hline
 \multicolumn{5}{l}{$P^1$ SLDG-time2-E}  \\
   5 &         3.16E-06 & -- &     3.01E-05 &      --         \\
   10 &        8.83E-06 &     1.48 &     5.06E-05 &     0.75 \\
   15 &        1.94E-05 &     1.95 &     9.47E-05 &     1.55 \\
   20 &       3.52E-05 &     2.06 &     1.64E-04 &     1.91  \\
   25 &       5.03E-05 &     1.60 &     2.30E-04 &     1.52 \\
   \hline
\multicolumn{5}{l}{$P^2$ SLDG-QC-time3-E }  \\
   5 &         1.26E-07 & -- &     1.23E-06 &   --         \\
   10 &        6.60E-07 &     2.39 &     5.59E-06 &     2.18\\
   15 &        2.20E-06 &     2.97 &     1.54E-05 &     2.50 \\
   20 &        5.26E-06 &     3.03 &     3.69E-05 &     3.04 \\
   25 &        9.61E-06 &     2.70 &     6.69E-05 &     2.67 \\
   \hline
\end{tabular}
\label{temporal_two1}
\end{table}

\begin{figure}[h!]
\centering
\includegraphics[height=65mm]{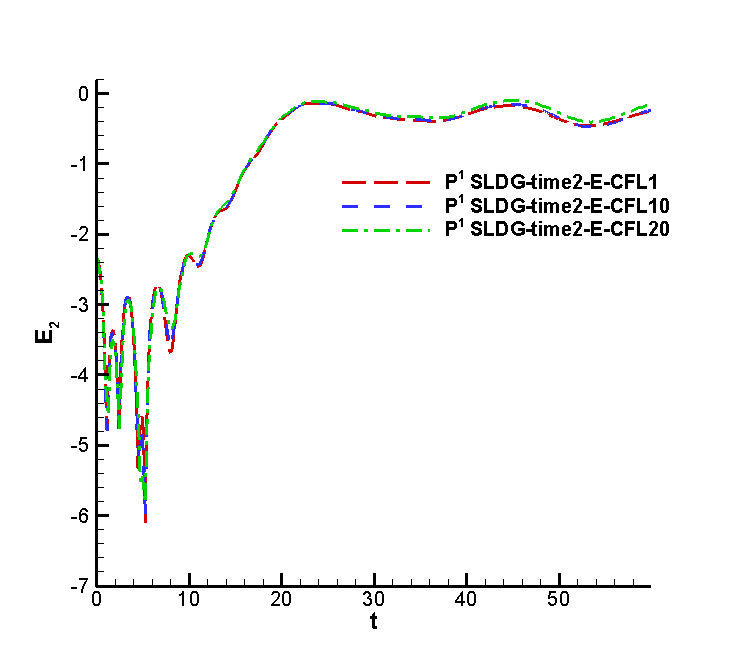}
\includegraphics[height=65mm]{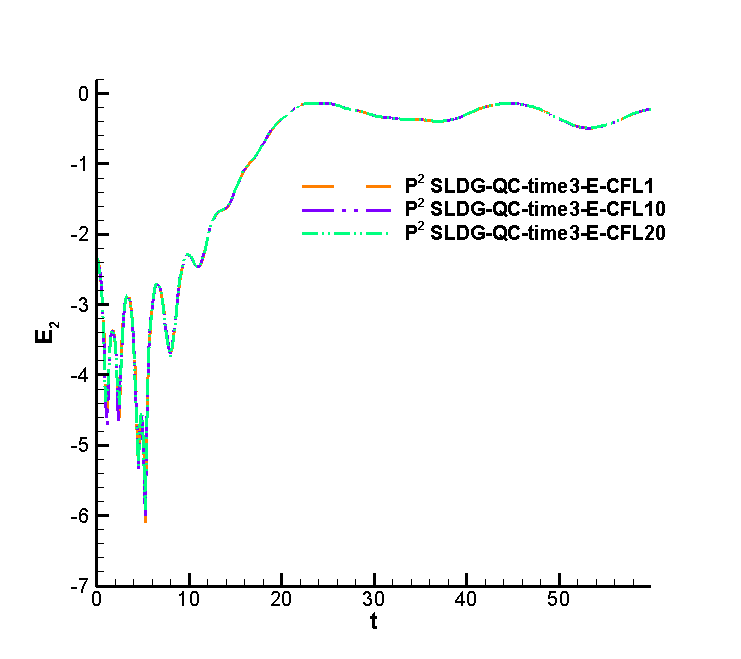}
\caption{Two-stream instability I. Time evolution of the electric field in $L^2$. }
\label{two1_e2}
\end{figure}

\begin{figure}[h!]
\centering
\includegraphics[height=65mm]{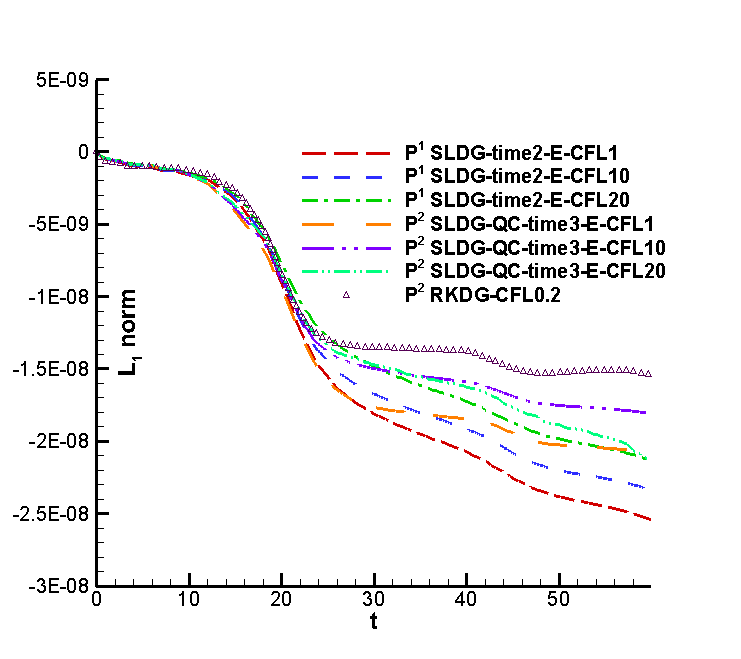}
\includegraphics[height=65mm]{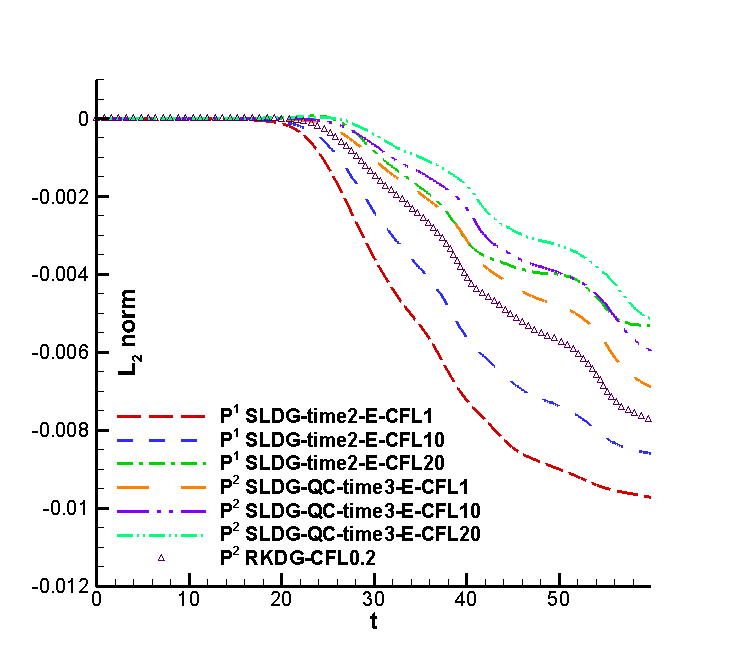}
\includegraphics[height=65mm]{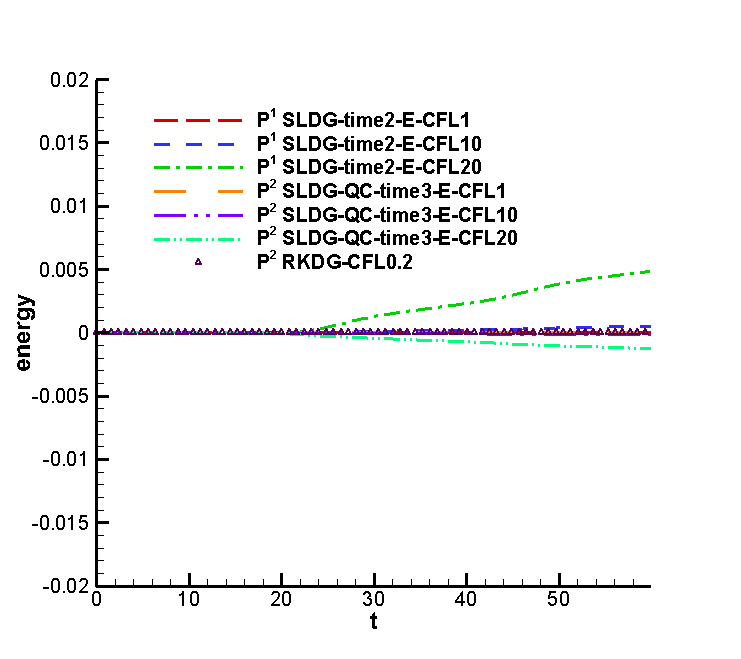}
\includegraphics[height=65mm]{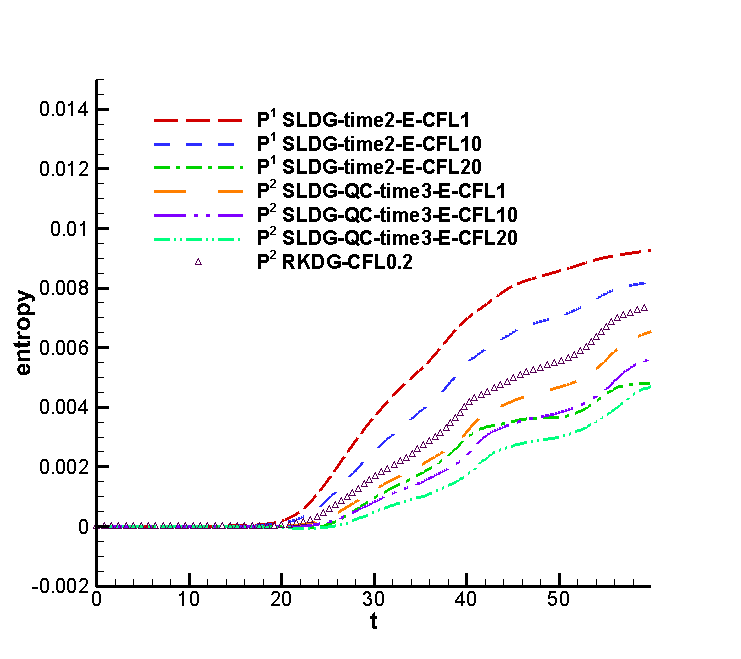}
\caption{Two-stream instability I. Time evolution of $L^1$ (upper left) and $L^2$ (upper right) norms of the solution as well as the discrete kinetic energy (lower left) and entropy (lower right).}
\label{two1norm}
\end{figure}

\begin{figure}[h!]
\centering
\subfigure[$P^2$ SLDG-QC-time3-E with $CFL=10$]{
\includegraphics[height=45mm]{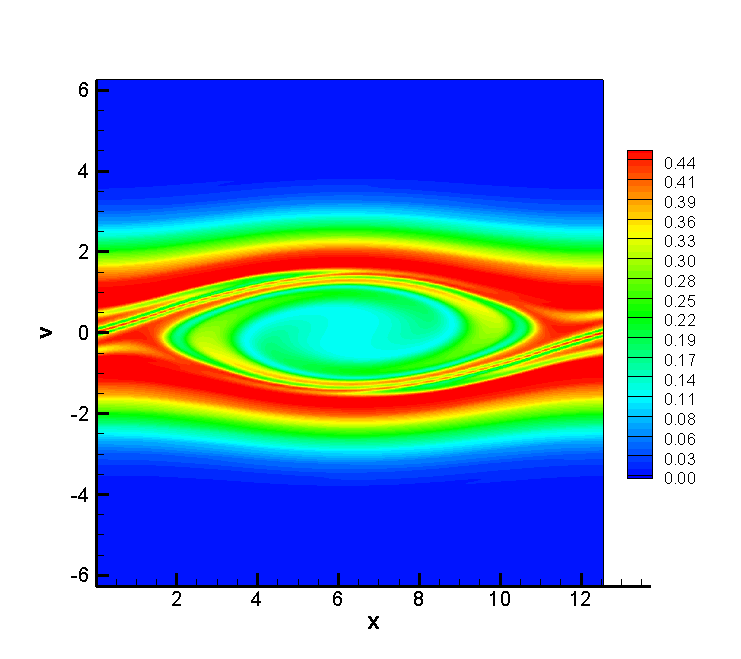} }
\subfigure[$P^2$ SLDG-QC-time3-E with $CFL=30$]{
\includegraphics[height=45mm]{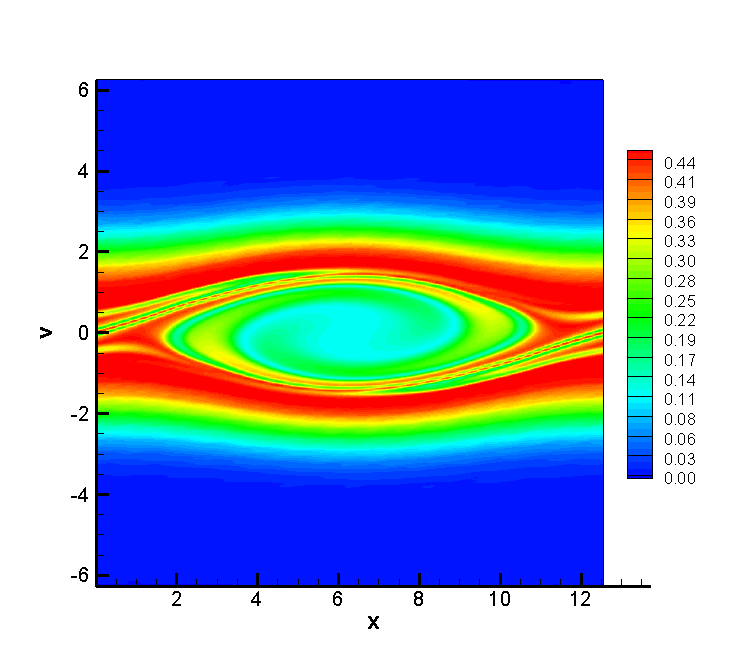} }
\subfigure[$P^2$ RKDG with $CFL=0.2$]{
\includegraphics[height=45mm]{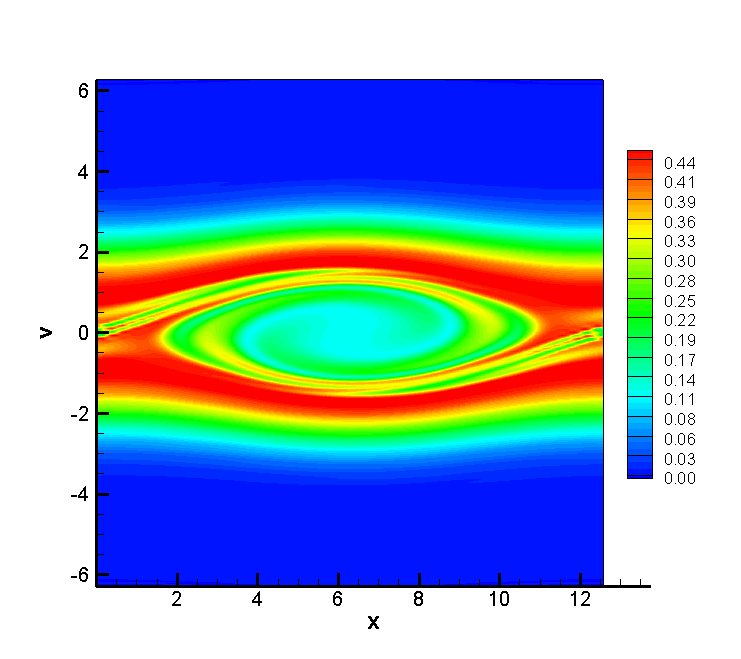} }

\caption{Two stream instability I at $T=40$.  Mesh: $160\times160$.
(a): $P^2$ SLDG-QC-time3-E with $CFL=10$. (b): $P^2$ SLDG-QC-time3-E with $CFL=30$.
(c): $P^2$ RKDG with $CFL=0.2$.}
\label{two1solution}
\end{figure}

\end{exa}

\begin{exa}
\emph{(Two Stream instability II.)}
Consider the symmetric two stream instability  \cite{umeda2008conservative,crouseilles2009conservative}, with the perturbed equilibrium as  the initial condition
\begin{equation*}
f(x,v,t=0) = \frac{1}{2v_t\sqrt{2\pi} } \left[ \exp\left( - \frac{(v-u)^2}{2v_{th}^2} \right)
+ \exp\left( - \frac{(v+u)^2}{2v_{th}^2} \right) \right] (1+0.05\cos(kx) ).
\end{equation*}
where $u=0.99$, $k=\frac{2}{13}$, and $v_t =0.3$.
We let $v_{\max} =2\pi$ and use a spatial mesh of $160\times 160$ elements.
We plot the time evolution of the electric field in the $L^2$ and $L^\infty$ norms (in semi-log scale) in Figure \ref{two2_e2}, which is benchmarked against the results reported in the literature. Time evolution of the relative derivation of the discrete $L^1$ norm, $L^2$ norm, energy and entropy in Figure \ref{two2norm}.
 Figure \ref{two2solution} shows the numerical solutions of phase space profiles computed by the $P^2$ SLDG-QC-time3-E method with $CFL=10,15$ and the $P^2$ RKDG method with $CFL=0.2$  at $T=40$.
 Decent numerical performance of the $P^2$ SLDG-QC-time3-E method with $CFL=10$ and the $P^2$ RKDG method with $CFL=0.2$ is observed.
 On the other hand, the $P^2$ SLDG-QC-time3-E method with a large $CFL=15$ is still stable and generates reasonable result, but some mild wiggles are observed.
Note that, for this example, distortion of approximate upstream cells (hence the break down of the code) is observed for the $P^2$ SLDG-QC-time3-E method with $CFL=20$. 

\begin{figure}[h!]
\centering
\includegraphics[height=65mm]{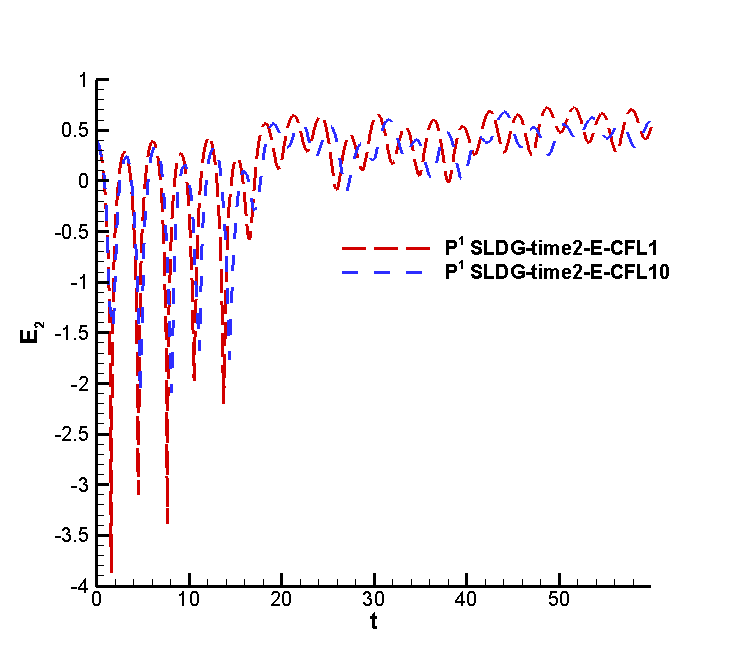}
\includegraphics[height=65mm]{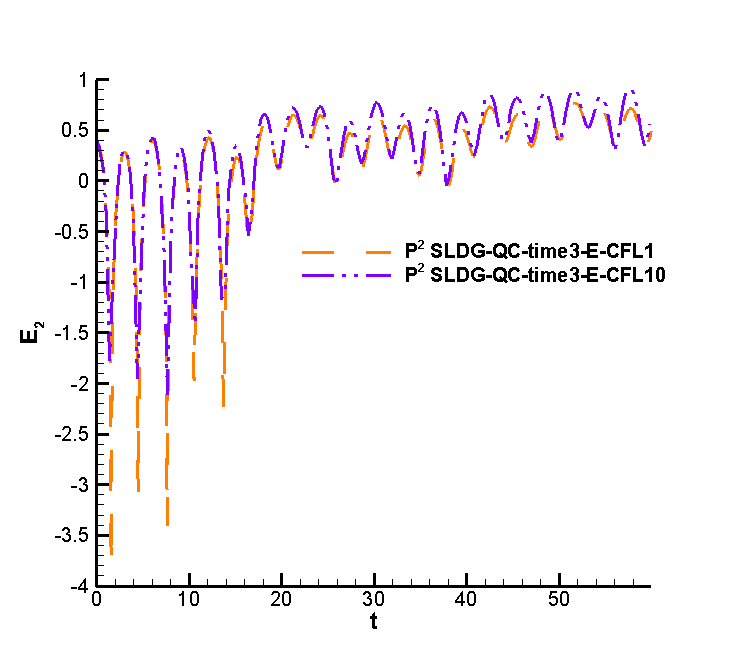}
\caption{Two stream instability II. Time evolution of the electric field in $L^2$. }
\label{two2_e2}
\end{figure}

\begin{figure}[h!]
\centering
\includegraphics[height=65mm]{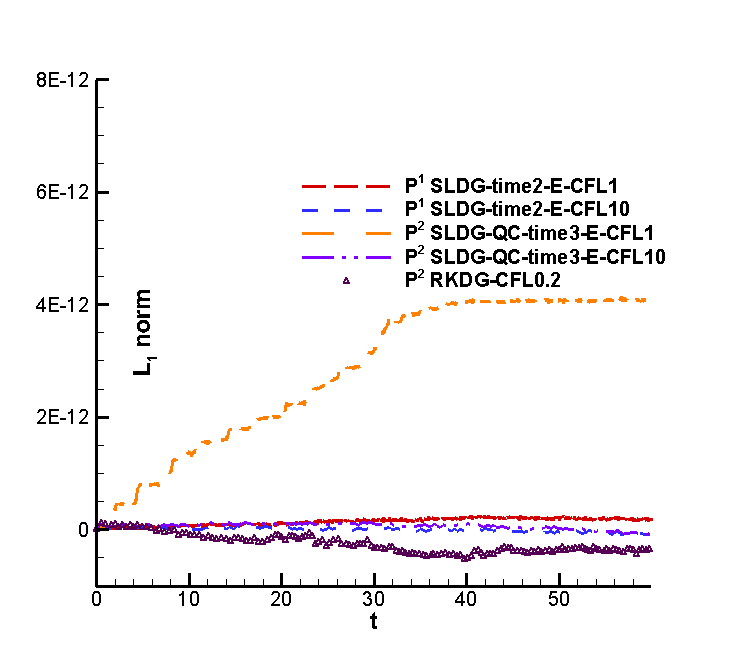}
\includegraphics[height=65mm]{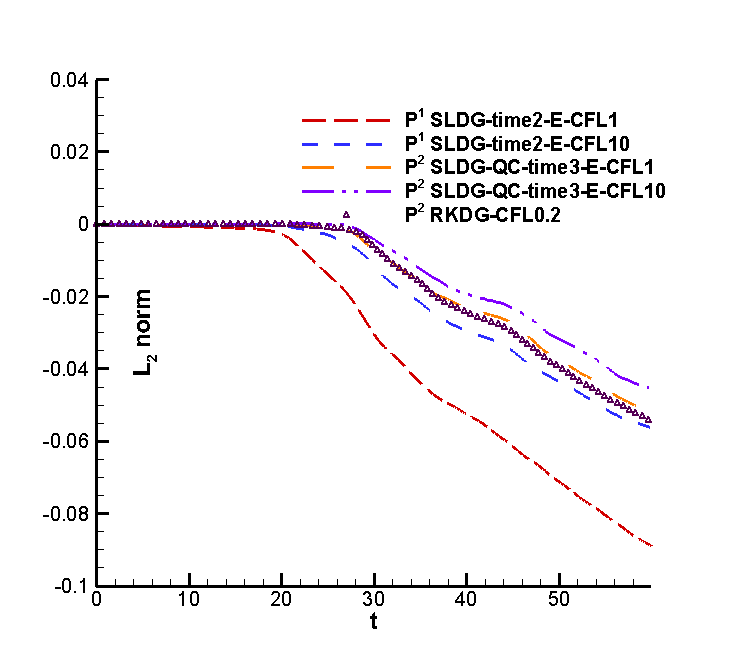}
\includegraphics[height=65mm]{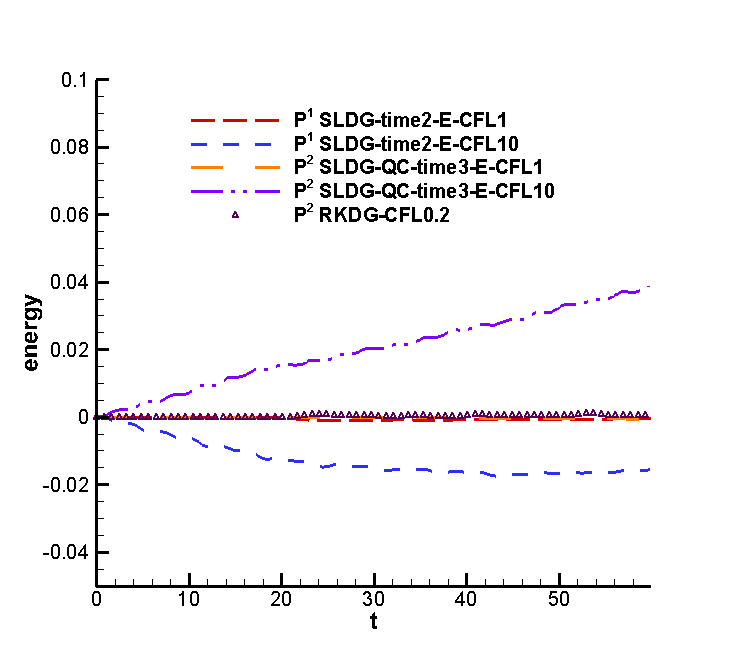}
\includegraphics[height=65mm]{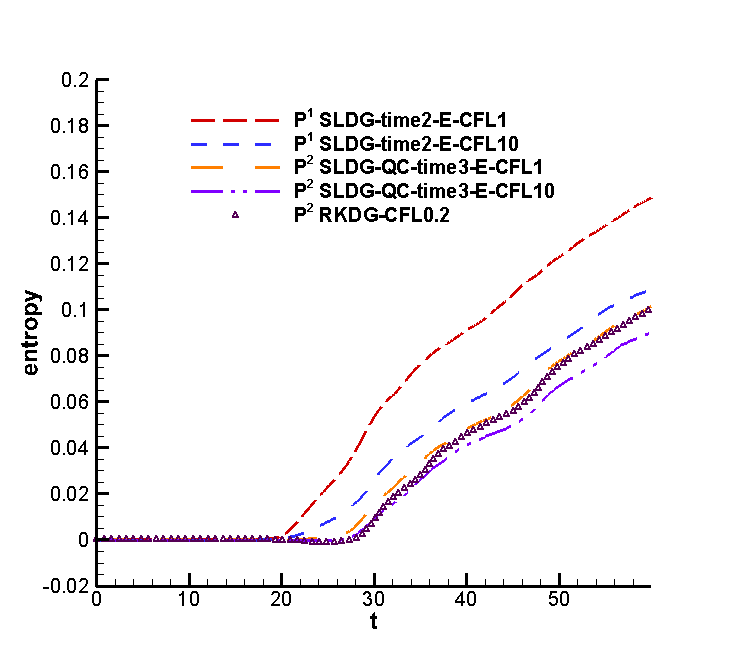}
\caption{Two stream instability II. Time evolution of $L^1$ (upper left) and $L^2$ (upper right) norms of the solution as well as the discrete kinetic energy (lower left) and entropy (lower right).}
\label{two2norm}
\end{figure}

\begin{figure}[h!]
\centering
\subfigure[$P^2$ SLDG-QC-time3-E with $CFL=10$]{
\includegraphics[height=45mm]{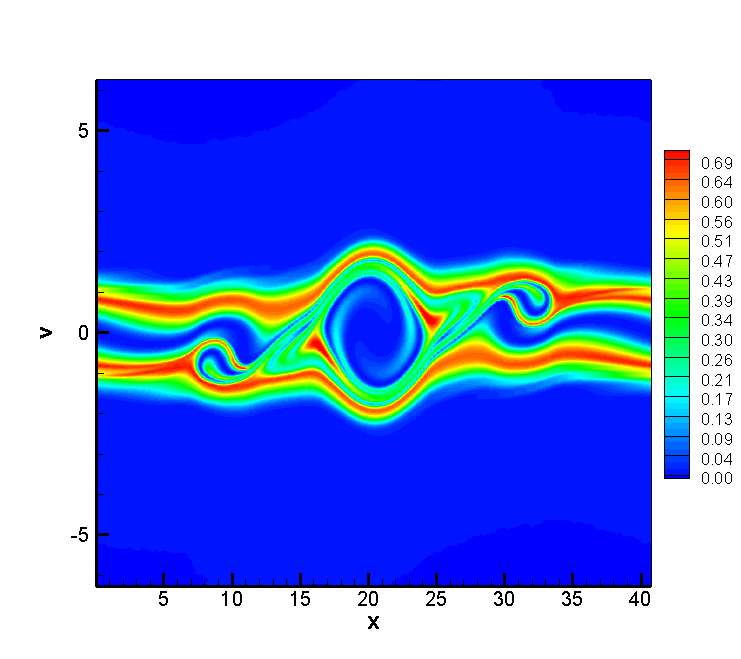} }
\subfigure[$P^2$ SLDG-QC-time3-E with $CFL=15$]{
\includegraphics[height=45mm]{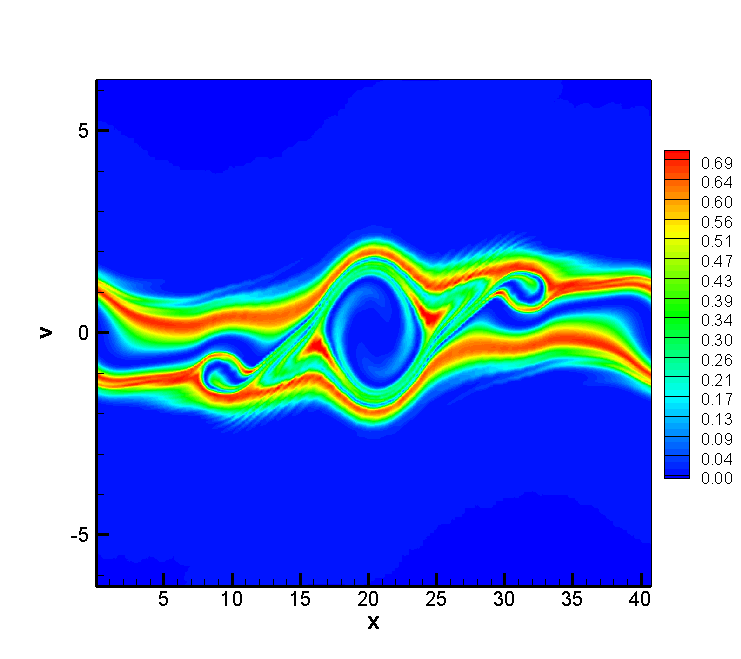} }
\subfigure[$P^2$ RKDG with $CFL=0.2$]{
\includegraphics[height=45mm]{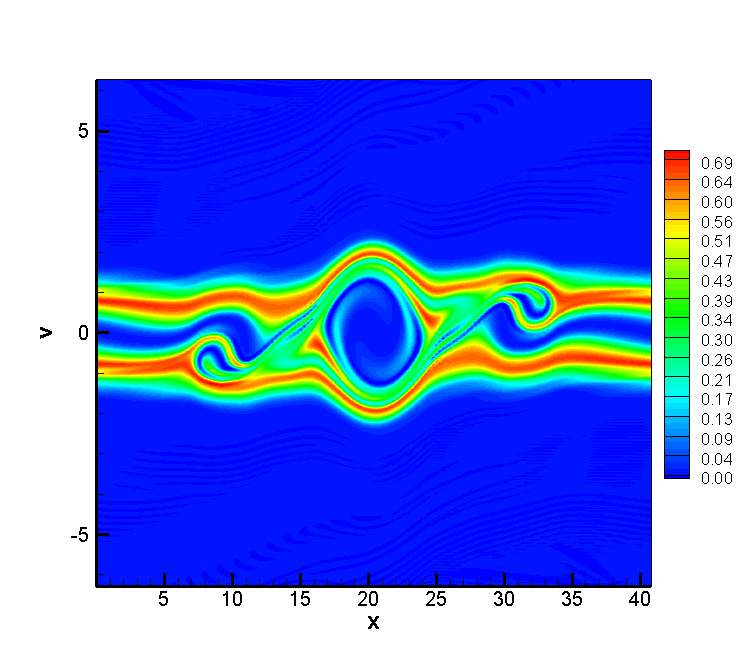} }
\caption{Two stream instability II at $T=40$. Mesh: $160\times160$.
(a): $P^2$ SLDG-QC-time3-E with $CFL=10$. (b): $P^2$ SLDG-QC-time3-E with $CFL=15$.
(c): $P^2$ RKDG with $CFL=0.2$. }
\label{two2solution}
\end{figure}
\end{exa}

\begin{exa}

\emph{(Bump-on-tail instability.)}
Consider an unstable bump-on-tail problem \cite{arber2002critical,xiong2014high}
with the initial distribution as
\begin{equation*}
f(x,v,t=0) = f_{BOT}(v) (1+0.04\cos(kx) ).
\end{equation*}
where the bump-on-tail distribution is
\begin{equation*}
f_{BOT} (v) = n_p \exp\left( - \frac{v^2}{2} \right)
+n_b \exp \left( - \frac{ (v-u)^2 }{2v_t^2} \right),
\end{equation*}
The parameters are chosen to be $n_p = \frac{9}{10\sqrt{2\pi} }$,  $n_b = \frac{2}{10\sqrt{2\pi} }$, $u=4.5$, $v_t =0.5$, $k=0.3$. The computational domain is $[0,\frac{20}{3}\pi] \times[ -13,13 ]$.
We show the time evolution of the electric field in $L^2$ norm (in semi-log scale) in Figure \ref{bump_e2}, and show the time evolution of the relative derivation of the discrete $L^1$ norm, $L^2$ norm, energy and entropy in Figure \ref{bump_norm}.
In Figure \ref{bump_solution}, we plot the numerical solutions of phase space profiles computed by the $P^2$ SLDG-QC-time3-E method with $CFL=10,30$ and the $P^2$ RKDG method with $CFL=0.2$  at $T=40$.
As in Example \ref{exa:two1}, the proposed SLDG schemes with $CFL$ as large as 30 are still able to generate very consistent results with that by the  $P^2$ RKDG method with $CFL=0.2$, leading to great computational savings.

\begin{figure}[h!]
\centering
\includegraphics[height=65mm]{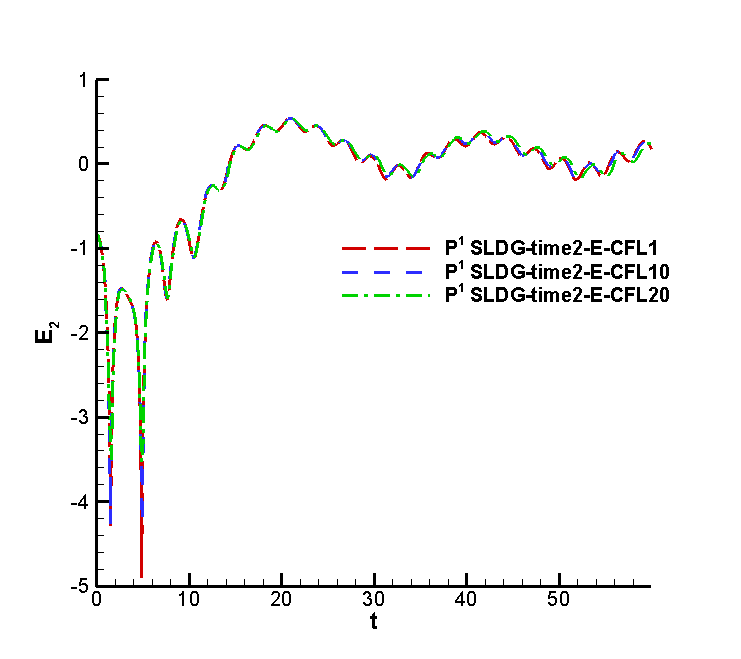}
\includegraphics[height=65mm]{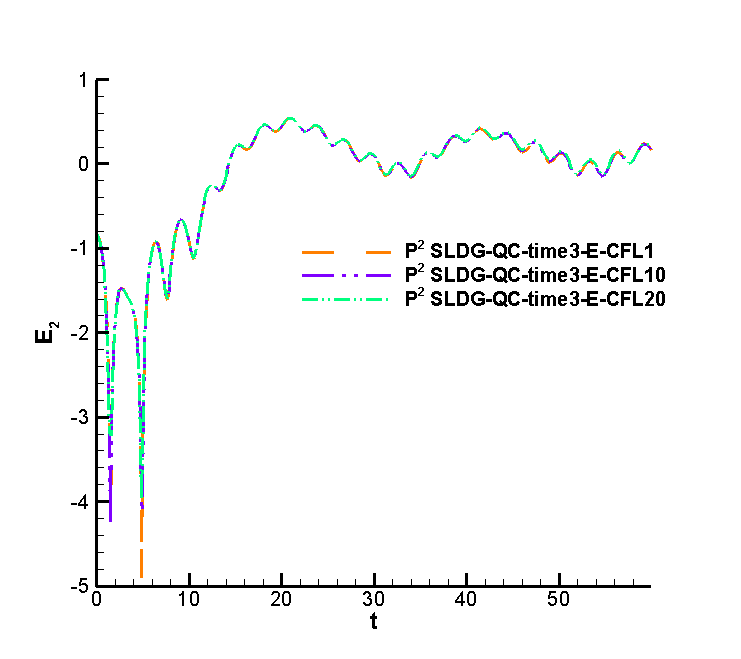}
\caption{Bump-on-tail instability. Time evolution of the electric field in $L^2$. }
\label{bump_e2}
\end{figure}

\begin{figure}[h!]
\centering
\includegraphics[height=65mm]{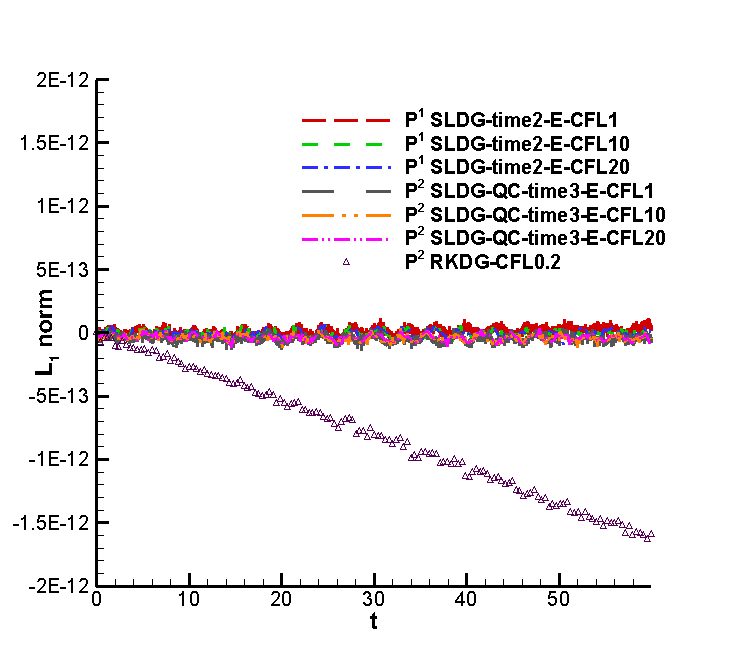}
\includegraphics[height=65mm]{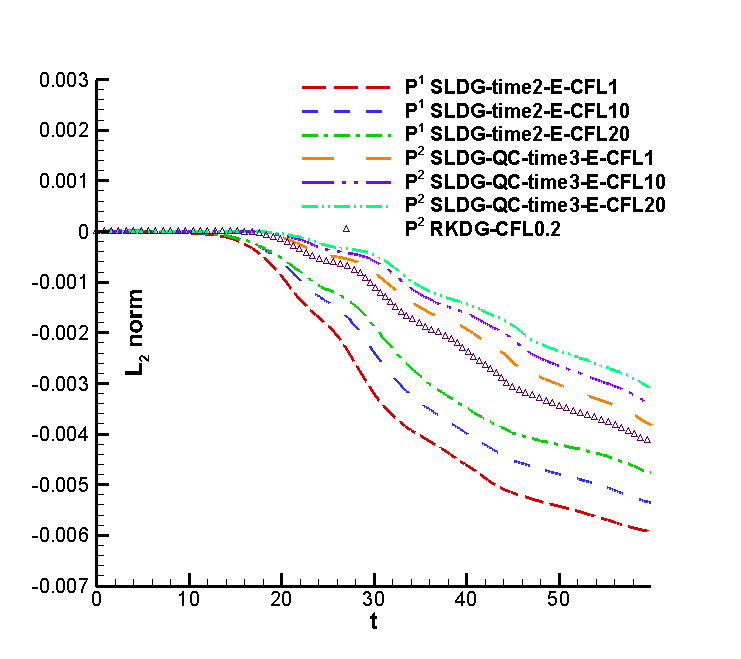}
\includegraphics[height=65mm]{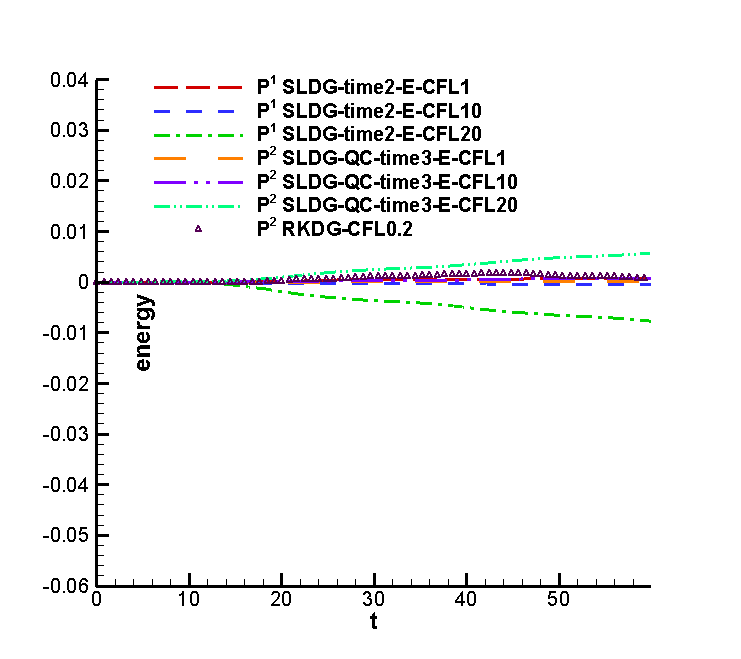}
\includegraphics[height=65mm]{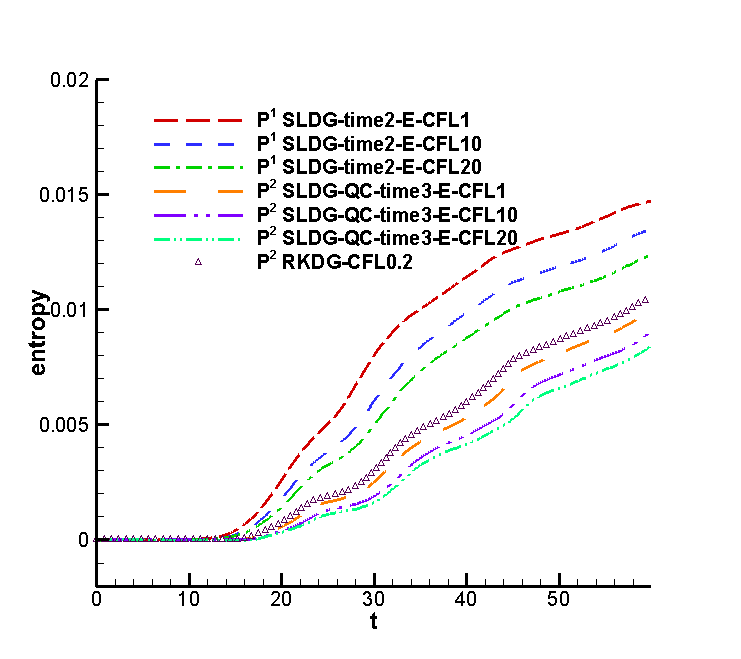}
\caption{Bump-on-tail instability. Time evolution of $L^1$ (upper left) and $L^2$ (upper right) norms of the solution as well as the discrete kinetic energy (lower left) and entropy (lower right).}
\label{bump_norm}
\end{figure}

\begin{figure}[h!]
\centering
\subfigure[$P^2$ SLDG-QC-time3-E with $CFL=10$]{
\includegraphics[height=45mm]{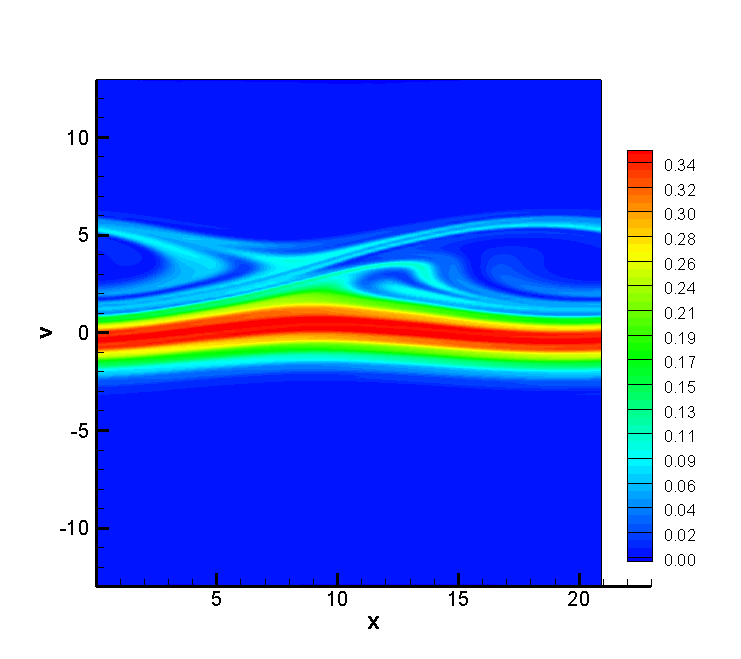} }
\subfigure[$P^2$ SLDG-QC-time3-E with $CFL=30$]{
\includegraphics[height=45mm]{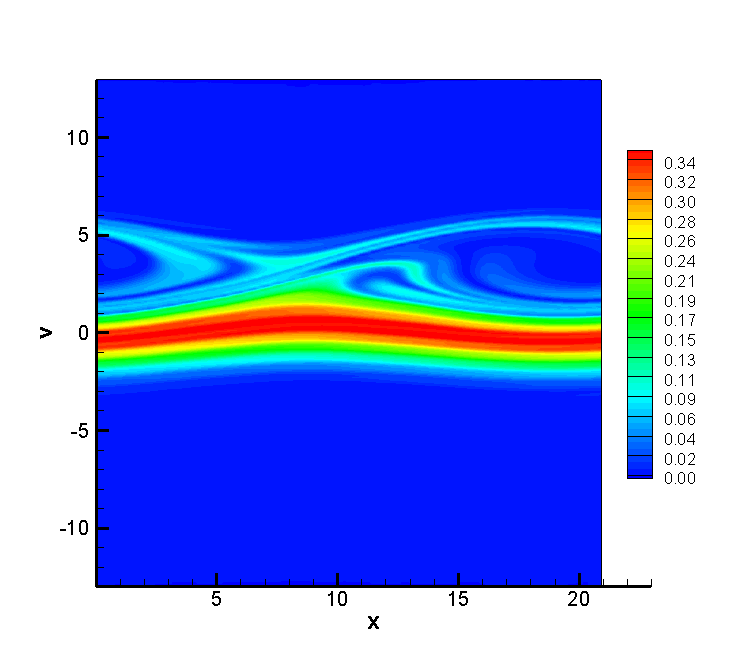} }
\subfigure[$P^2$ RKDG with $CFL=0.2$]{
\includegraphics[height=45mm]{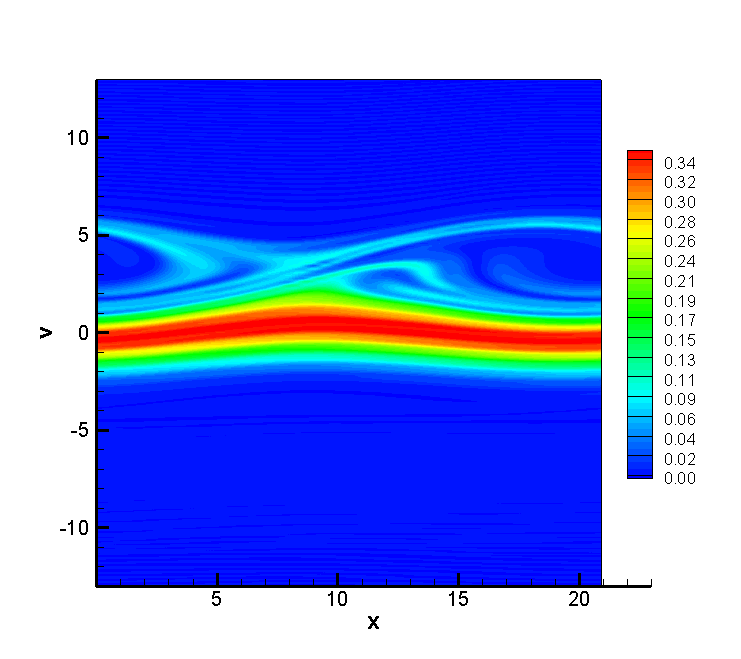} }
\caption{Bump-on-tail instability.  Mesh: $160\times160$.
(a): $P^2$ SLDG-QC-time3-E with $CFL=10$. (b): $P^2$ SLDG-QC-time3-E with $CFL=30$.
(c): $P^2$ RKDG with $CFL=0.2$. }
\label{bump_solution}
\end{figure}

\end{exa}

%% file: conclusion.tex
\section{Conclusion}
 A high order SLDG method was proposed for solving the VP system. The two key ingredients of the proposed scheme are (1) a high order non-splitting conservative SLDG transport scheme and (2) a high order characteristics tracking approach for the VP system. The proposed method is locally mass conservative, highly accurate,  free of splitting error and allows for extra large time stepping size. To the best of the authors' knowledge, this is the first SLDG scheme that is able to attain all these desired properties. The numerical performance of the method is promising. 
We compare the CPU time of an efficient implementation of the proposed method, with that of the Eulerian RK DG method, in many benchmark VP test problems. Tremendous computational savings are observed, without compromising effectiveness of the method. We will perform rigorous error estimate, investigate superconvergence properties, and generalize the scheme to handle non-trivial boundary conditions, diffusion terms, source terms, etc. in our future work. 
The extension to the VP system in higher dimensions will be explored as well.

%% file: appendix.tex
\section{Search intersection points between quadratic-curved sides of an upstream cell and background grid lines}
\label{append:a}
The algorithm for searching intersection points between an upstream cell and the grid line $x= x_i$ is described in the following. The procedure for searching intersection points between the upstream cell and $v= v_j$ is pretty similar, thus omitted for brevity. 
    We can find intersection points of the quadratic curve determined by points $(x^\star_1,v^\star_1)$, $(x^\star_2,v^\star_2)$, $(x^\star_3,v^\star_3)$  and grid line $x=x_i$ by solving the following equation,
         \begin{equation}
     \begin{cases}
     x_i =   \frac{x^\star_3 -x^\star_1 }{2} \xi + \frac{v^\star_3 -v^\star_1 }{2} \eta
       + \frac{ x^\star_3 +x^\star_1 }{2}  , \\
      \eta = \frac{\eta_2}{ \xi_2^2-1 }(  \xi^2-1  ).
     \end{cases}
     \label{x_i}
     \end{equation}

     \begin{description}
       \item[Case 1]
       If $|x^\star_3- x^\star_1| \leq |v^\star_3 - v^\star_1| $, we have the following equation derived from \eqref{x_i},
       \begin{equation}
         A\xi^2 +  B\xi +C
       =0.
       \label{xi_quadratic}
       \end{equation}
    where
       \begin{align*}
       A&= \frac{\eta_2}{ \xi_2^2 - 1 }, \quad
       B= \frac{ x^\star_3 -x^\star_1 }{ v^\star_3 -v^\star_1 }, \quad
       C= - \frac{ \left(x_i -\frac{x^\star_1+x^\star_3}{2}\right) }{  \frac{v^\star_3-v^\star_1}{2} } -\frac{\eta_2 }{ \xi_2^2 -1 }.
       \end{align*}

       We can solve \eqref{xi_quadratic} as follows:
       
  \bigskip     
\fbox{
\begin{minipage}[htb]{0.9\linewidth}
         \begin{algorithmic}
  \IF{ $A\ge \varepsilon$, where $\varepsilon=10^{-13}$ }
      \STATE
      Let
      $\Delta = B^2 -4AC$.

      If $\Delta<0$, there is no solution.

      If $\Delta=0$, there is only one solution,
      $\xi_1 = - \frac{B}{2A}.$

%

%

       If $\Delta>0$, there are two solutions,
$$
      \xi_{1} =\frac{ 2C }{ -B - \gamma \sqrt{ \Delta } }, \quad
      \xi_2 =\frac{ -B - \gamma \sqrt{ \Delta } }{2A},
$$
      where $\gamma=1$ if $B\geq 0$, and $\gamma=-1$ otherwise.
 \ELSIF{ $A<\varepsilon$ and $B \ge \varepsilon$ }
  \STATE there is only one solution, $\xi_1 = - \frac{C}{B}.$
  \ELSE
  \STATE We retreat this case as no intersection points.


  \ENDIF
  
       If $\xi\in[-1,1]$, the solution $(\xi,\eta)$ is identified as an intersection point.
\end{algorithmic}
\end{minipage}
}

       \item[Case 2]
 If $|x^\star_3- x^\star_1| > |v^\star_3 - v^\star_1| $, we have the following equation derived from \eqref{x_i},
       \begin{equation*}
       A\eta^2 +B\eta +C =0,
       \end{equation*}
       where
       \begin{align*}
       A&= \frac{\eta_2}{\xi_2^2 - 1} \left( \frac{v^\star_3-v^\star_1}{x^\star_3-x^\star_1} \right)^2,\\
       B&=   -1- \frac{ 4\left(x_i-\frac{x^\star_1+x^\star_3}{2} \right) }{ x^\star_2 -x^\star_1 } \frac{v^\star_3 -v^\star_1 }{x^\star_3 -x^\star_1} \frac{ \eta_2 }{\xi_2^2-1},\\
       C&= \frac{\eta_2}{\xi_2^2 - 1}
       \left(  \frac{ 4\left( x_i - \frac{x^\star_1+x^\star_3}{2} \right)^2 }{ (x^\star_3 -x^\star_1)^2 }  -1 \right).
       \end{align*}

       Similar to Case 1.1, the solution $\eta$ can be solved. And then,
       \begin{equation}
       \xi = \frac{ 2x_i - x^\star_3 - x^\star_1 }{ x^\star_3 - x^\star_1 } - \frac{ v^\star_3 -v^\star_1 }{x^\star_3 -x^\star_1}\eta.
       \end{equation}

       If $\xi\in[-1,1]$, the solution $(\xi,\eta)$ is identified as an intersection point.

     \end{description}

%% file: main.bbl
\begin{thebibliography}{10}

\bibitem{arber2002critical}
T.~Arber and R.~Vann.
\newblock {A critical comparison of Eulerian-grid-based Vlasov solvers}.
\newblock {\em Journal of computational physics}, 180(1):339--357, 2002.

\bibitem{arnold2002unified}
D.~Arnold, F.~Brezzi, B.~Cockburn, and L.~Marini.
\newblock {Unified analysis of discontinuous Galerkin methods for elliptic
  problems}.
\newblock {\em SIAM Journal on Numerical Analysis}, 39(5):1749--1779, 2002.

\bibitem{birdsall2005plasma}
C.~Birdsall and A.~Langdon.
\newblock {\em Plasma physics via computer simulaition}.
\newblock CRC Press, 2005.

\bibitem{cai2016high}
X.~Cai, W.~Guo, and J.-M. Qiu.
\newblock {A high order conservative semi-Lagrangian discontinuous Galerkin
  method for two-dimensional transport simulations}.
\newblock {\em Journal of Scientific Computing}, accepted, 2017.

\bibitem{cai2016conservative}
X.~Cai, J.~Qiu, and J.-M. Qiu.
\newblock {A conservative semi-Lagrangian HWENO method for the Vlasov
  equation}.
\newblock {\em Journal of Computational Physics}, 323:95--114, 2016.

\bibitem{carrillo2007nonoscillatory}
J.~Carrillo and F.~Vecil.
\newblock {Nonoscillatory interpolation methods applied to Vlasov-based
  models}.
\newblock {\em SIAM Journal on Scientific Computing}, 29(3):1179--1206, 2007.

\bibitem{casas2016high}
F.~Casas, N.~Crouseilles, E.~Faou, and M.~Mehrenberger.
\newblock {High-order Hamiltonian splitting for the Vlasov--Poisson equations}.
\newblock {\em Numerische Mathematik}, 135(3):769--801, 2017.

\bibitem{castillo2000priori}
P.~Castillo, B.~Cockburn, I.~Perugia, and D.~Sch{\"o}tzau.
\newblock {An a priori error analysis of the local discontinuous Galerkin
  method for elliptic problems}.
\newblock {\em SIAM Journal on Numerical Analysis}, 38(5):1676--1706, 2000.

\bibitem{cheng}
C.~Cheng and G.~Knorr.
\newblock {The integration of the Vlasov equation in configuration space}.
\newblock {\em Journal of Computational Physics}, 22(3):330--351, 1976.

\bibitem{cheng2014energy}
Y.~Cheng, A.~Christlieb, and X.~Zhong.
\newblock {Energy-conserving discontinuous Galerkin methods for the
  {V}lasov--{A}mp\`{e}re system}.
\newblock {\em Journal of Computational Physics}, 256:630--655, 2014.

\bibitem{cheng2014energyvm}
Y.~Cheng, A.~Christlieb, and X.~Zhong.
\newblock {Energy-conserving discontinuous Galerkin methods for the
  Vlasov-Maxwell system}.
\newblock {\em Journal of Computational Physics}, 279:145--173, 2014.

\bibitem{cheng2014discontinuous}
Y.~Cheng, I.~Gamba, F.~Li, and P.~Morrison.
\newblock {Discontinuous Galerkin methods for the Vlasov--Maxwell equations}.
\newblock {\em SIAM Journal on Numerical Analysis}, 52(2):1017--1049, 2014.

\bibitem{cheng2012study}
Y.~Cheng, I.~Gamba, and P.~Morrison.
\newblock Study of conservation and recurrence of {Runge--Kutta} discontinuous
  {G}alerkin schemes for {V}lasov--{P}oisson systems.
\newblock {\em Journal of Scientific Computing}, 56(2):319--349, 2013.

\bibitem{christlieb2014high}
A.~Christlieb, W.~Guo, M.~Morton, and J.-M. Qiu.
\newblock {A high order time splitting method based on integral deferred
  correction for semi-Lagrangian Vlasov simulations}.
\newblock {\em Journal of Computational Physics}, 267:7--27, 2014.

\bibitem{cockburn2000discontinuous}
B.~Cockburn, G.~Karniadakis, and C.-W. Shu.
\newblock {\em {Discontinuous Galerkin Methods, Theory, Computation and
  Applications}}.
\newblock DOE/ER/25372-Final, Brown University, Providence, RI (US), 2000.

\bibitem{cockburn1998local}
B.~Cockburn and C.-W. Shu.
\newblock {The local discontinuous Galerkin method for time-dependent
  convection-diffusion systems}.
\newblock {\em SIAM Journal on Numerical Analysis}, 35(6):2440--2463, 1998.

\bibitem{crouseilles2009conservative}
N.~Crouseilles, M.~Mehrenberger, and E.~Sonnendr{\"u}cker.
\newblock {Conservative semi-Lagrangian schemes for Vlasov equations}.
\newblock {\em Journal of Computational Physics}, 229(6):1927--1953, 2010.

\bibitem{filbet2003comparison}
F.~Filbet and E.~Sonnendrucker.
\newblock {Comparison of Eulerian Vlasov solvers}.
\newblock {\em Computer Physics Communications}, 150(3):247--266, 2003.

\bibitem{filbet2001conservative}
F.~Filbet, E.~Sonnendr{\"u}cker, and P.~Bertrand.
\newblock {Conservative numerical schemes for the Vlasov equation}.
\newblock {\em Journal of Computational Physics}, 172(1):166--187, 2001.

\bibitem{FilbetSB}
F.~Filbet, E.~Sonnendr{\"u}cker, and P.~Bertrand.
\newblock Conservative numerical schemes for the {V}lasov equation.
\newblock {\em Journal of Computational Physics}, 172(1):166--187, 2001.

\bibitem{gucclu2014arbitrarily}
Y.~G{\"u}{\c{c}}l{\"u}, A.~J. Christlieb, and W.~N. Hitchon.
\newblock {Arbitrarily high order Convected Scheme solution of the
  Vlasov--Poisson system}.
\newblock {\em Journal of Computational Physics}, 270:711--752, 2014.

\bibitem{Guo2013discontinuous}
W.~Guo, R.~Nair, and J.-M. Qiu.
\newblock A conservative semi-{L}agrangian discontinuous {G}alerkin scheme on
  the cubed-sphere.
\newblock {\em Monthly Weather Review}, 142(1):457--475, 2013.

\bibitem{guo2015efficient}
W.~Guo, R.~Nair, and X.~Zhong.
\newblock {An efficient WENO limiter for discontinuous Galerkin transport
  scheme on the cubed sphere}.
\newblock {\em International Journal for Numerical Methods in Fluids},
  81:3--21, 2015.

\bibitem{guo2013hybrid}
W.~Guo and J.-M. Qiu.
\newblock {Hybrid semi-Lagrangian finite element-finite difference methods for
  the Vlasov equation}.
\newblock {\em Journal of Computational Physics}, 234:108--132, 2013.

\bibitem{heath2012discontinuous}
R.~Heath, I.~Gamba, P.~Morrison, and C.~Michler.
\newblock {A discontinuous Galerkin method for the Vlasov--Poisson system}.
\newblock {\em Journal of Computational Physics}, 231(4):1140--1174, 2012.

\bibitem{lauritzen2010conservative}
P.~Lauritzen, R.~Nair, and P.~Ullrich.
\newblock {A conservative semi-Lagrangian multi-tracer transport scheme (CSLAM)
  on the cubed-sphere grid}.
\newblock {\em Journal of Computational Physics}, 229(5):1401--1424, 2010.

\bibitem{nakamura1999cubic}
T.~Nakamura and T.~Yabe.
\newblock Cubic interpolated propagation scheme for solving the
  hyper-dimensional {V}lasov--{P}oisson equation in phase space.
\newblock {\em Computer Physics Communications}, 120(2):122--154, 1999.

\bibitem{Qiu_Christlieb}
J.-M. Qiu and A.~Christlieb.
\newblock {A Conservative high order semi-Lagrangian WENO method for the Vlasov
  Equation}.
\newblock {\em Journal of Computational Physics}, 229:1130--1149, 2010.

\bibitem{qiu2017high}
J.-M. Qiu and G.~Russo.
\newblock {A High Order Multi-Dimensional Characteristic Tracing Strategy for
  the Vlasov--Poisson System}.
\newblock {\em Journal of Scientific Computing}, 71(1):414--434, 2017.

\bibitem{qiu2011conservative}
J.-M. Qiu and C.-W. Shu.
\newblock {Conservative semi-Lagrangian finite difference WENO formulations
  with applications to the Vlasov equation}.
\newblock {\em Communications in Computational Physics}, 10(4):979, 2011.

\bibitem{qiu2011positivity}
J.-M. Qiu and C.-W. Shu.
\newblock {Positivity preserving semi-Lagrangian discontinuous Galerkin
  formulation: Theoretical analysis and application to the Vlasov--Poisson
  system}.
\newblock {\em Journal of Computational Physics}, 230(23):8386--8409, 2011.

\bibitem{restelli2006semi}
M.~Restelli, L.~Bonaventura, and R.~Sacco.
\newblock {A semi-Lagrangian discontinuous Galerkin method for scalar advection
  by incompressible flows}.
\newblock {\em Journal of Computational Physics}, 216(1):195--215, 2006.

\bibitem{rossmanith2011positivity}
J.~A. Rossmanith and D.~C. Seal.
\newblock A positivity-preserving high-order semi-{L}agrangian discontinuous
  {G}alerkin scheme for the {V}lasov--{P}oisson equations.
\newblock {\em Journal of Computational Physics}, 230(16):6203--6232, 2011.

\bibitem{sonnendruecker}
E.~Sonnendruecker, J.~Roche, P.~Bertrand, and A.~Ghizzo.
\newblock {The semi-Lagrangian method for the numerical resolution of the
  Vlasov equation}.
\newblock {\em Journal of Computational Physics}, 149(2):201--220, 1999.

\bibitem{umeda2008conservative}
T.~Umeda.
\newblock {A conservative and non-oscillatory scheme for Vlasov code
  simulations}.
\newblock {\em Earth, planets and space}, 60(7):773--779, 2008.

\bibitem{xiong2014high}
T.~Xiong, J.-M. Qiu, Z.~Xu, and A.~Christlieb.
\newblock {High order maximum principle preserving semi-Lagrangian finite
  difference WENO schemes for the Vlasov equation}.
\newblock {\em Journal of Computational Physics}, 273:618--639, 2014.

\bibitem{xiong2016conservative}
T.~Xiong, G.~Russo, and J.-M. Qiu.
\newblock {Conservative multi-dimensional semi-Lagrangian finite difference
  scheme: stability and applications to the kinetic and fluid simulations}.
\newblock {\em arXiv preprint arXiv:1607.07409}, 2016.

\bibitem{zhangshu2010}
X.~Zhang and C.-W. Shu.
\newblock {On maximum-principle-satisfying high order schemes for scalar
  conservation laws}.
\newblock {\em Journal of Computational Physics}, 229:3091--3120, 2010.

\end{thebibliography}
